\documentclass[a4paper,twoside,11pt,reqno]{amsart}
\newcommand{\Ueberschrift}{Weil--Ch\^atelet divisible elements in Tate--Shafarevich groups II: \\ On a question of Cassels}
\newcommand{\Kurztitel}{On Weil--Ch\^atelet divisible elements in Tate--Shafarevich groups}
\usepackage[headings]{fullpage}
\usepackage{amsmath} \usepackage{amssymb} \usepackage{amsopn}
\usepackage{amsthm} \usepackage{amsfonts} \usepackage{mathrsfs} 
\usepackage[dvips]{graphicx} \usepackage[usenames]{color}
\usepackage[all]{xy}
\usepackage[OT2, T1]{fontenc}
\usepackage[latin1]{inputenc}  
\usepackage[pdfpagelabels, pdftex]{hyperref}

\usepackage{dsfont}

\hypersetup{
  pdftitle={},
  pdfauthor={},
  pdfsubject={},
  pdfkeywords={},
  colorlinks=true,
  breaklinks=true,
  bookmarksopen=true,
  bookmarksnumbered=true,
  pdfpagemode=UseOutlines,
  plainpages=false}

\DeclareMathOperator{\rH}{H}
\DeclareMathOperator{\rI}{I}

\DeclareMathOperator{\rT}{T}

\DeclareMathOperator{\rc}{c}

\newcommand{\bF}{{\mathbb F}}

\newcommand{\bL}{{\mathbb L}}

\newcommand{\bN}{{\mathbb N}}

\newcommand{\bQ}{{\mathbb Q}}

\newcommand{\bZ}{{\mathbb Z}}

\newcommand{\fo}{{\mathfrak o}}

\DeclareSymbolFont{cyrletters}{OT2}{wncyr}{m}{n}
\DeclareMathSymbol{\Sha}{\mathalpha}{cyrletters}{"58}

\newcommand{\one}{\mathbf{1}}
\newcommand{\surj}{\twoheadrightarrow} 
\newcommand{\inj}{\hookrightarrow}

\newcommand{\ev}{{\rm ev}}

\DeclareMathOperator{\Hom}{Hom}

\DeclareMathOperator{\Aut}{Aut}
\DeclareMathOperator{\End}{End}

\DeclareMathOperator{\im}{im}

\DeclareMathOperator{\GL}{GL}
\DeclareMathOperator{\PGL}{PGL}
\DeclareMathOperator{\SL}{SL}

\newcommand{\tr}{{\rm tr}} 

\newcommand{\matzz}[4]{\left(
\begin{array}{cc} #1 & #2 \\ #3 & #4 \end{array} \right)}

\DeclareMathOperator{\Spec}{Spec}

\DeclareMathOperator{\Div}{Div}
\DeclareMathOperator{\divisor}{div}
\DeclareMathOperator{\Pic}{Pic}

\DeclareMathOperator{\Frob}{Frob}

\DeclareMathOperator{\res}{res}

\DeclareMathOperator{\infl}{inf}
\DeclareMathOperator{\ind}{ind}

\DeclareMathOperator{\Gal}{Gal}

\newcommand{\ph}{\varphi}

\newcommand{\alg}{{\rm alg}}
\newcommand{\nr}{{\rm nr}}

\newcommand{\ab}{{\rm ab}}

\newcommand{\Sel}{{\rm Sel}} 
\newcommand{\kum}{{\rm kum}}

\newcommand{\bruch}[2]{\genfrac{}{}{0.5pt}{}{#1}{#2}}

\newcommand{\ov}[1]{\mbox{${\overline{#1}}$}} 

\newtheorem{thm}{Theorem}

\newtheorem{prop}[thm]{Proposition}
\newtheorem{lem}[thm]{Lemma}

\newtheorem{cor}[thm]{Corollary}

\newtheorem*{thmA}{Theorem A}
\newtheorem*{thmB}{Theorem B}
\newtheorem*{thmC}{Theorem C}
\newtheorem*{thmD}{Theorem D}
\newtheorem*{thmD'}{Theorem D'}

\theoremstyle{definition}

\theoremstyle{remark}
\newtheorem{rmk}[thm]{Remark}

\newtheorem{ex}[thm]{Example}

\newenvironment{pro}[1][Proof]{{\it{#1:}} }{\hfill $\square$}
\newenvironment{pro*}[1][Proof]{{\it{#1:}} }{}

\numberwithin{equation}{section}

\begin{document}

\title[\Kurztitel]{\Ueberschrift} 
\author{Mirela \c{C}iperiani}
\address{Mirela \c{C}iperiani, Department of Mathematics, the University of Texas at Austin, 1 University Station, C1200 
Austin, Texas 78712, USA}
\email{mirela@math.utexas.edu}
\urladdr{http://www.ma.utexas.edu/users/mirela/}

\author{Jakob Stix}
\address{Jakob Stix, Mathematisches  Institut, Universit\"at Heidelberg, Im Neuenheimer Feld 288, 69120 Heidelberg, Germany}
\email{stix@mathi.uni-heidelberg.de}
\urladdr{http://www.mathi.uni-heidelberg.de/~stix/}
\thanks{The authors acknowledge the hospitality and support provided by MATCH and the Newton Institute. \\
\indent The first author was partially supported by NSF Grant DMS-07-58362 and by NSA Grant H98230-12-1-0208.}

\date{February 19, 2013}

\maketitle

\begin{quotation}
\noindent \small {\bf Abstract} --- For an abelian variety $A$ over a number field $k$ we discuss the divisibility in $\rH^1(k,A)$ of elements of the subgroup $\Sha(A/k)$. The results are most complete for elliptic curves over $\bQ$. 
\end{quotation}

%%%%%%%%%%%%%%%%%%%%%%%%%%%%%%%%%%%%%%%%%%%
%%%%%%%%%%%%%%%%%%%%%%%%%%%%%%%%%%%%%%%%%%%

\setcounter{tocdepth}{1} {\scriptsize \tableofcontents}

%%%%%%%%%%%%%%%%%%%%%%%%%%%%%%%%%%%%%%%%%%%
%%%%%%%%%%%%%%%%%%%%%%%%%%%%%%%%%%%%%%%%%%%

\section{Introduction}

%-----------------------------------------------------------------------------------------------------------------------

\noindent Let $A$ be an abelian variety over an algebraic number field $k$ with absolute Galois group $\Gal_k$. We aim to determine whether elements of the Tate--Shafarevich group
\[
\Sha(A/k)
\]
can become divisible in the  Weil--Ch\^atelet group $\rH^1(k,A)= \rH^1(k,A(k^\alg))$ with $k^\alg$ the algebraic closure of $k$, i.e., lie in the subgroup of divisible elements
\[
\divisor(\rH^1(k,A)) = \bigcap_{n\in \bN} n\rH^1(k,A).
\]
This question was initially asked by Cassels in the case of elliptic curves (see \cite{cassels:III}  Problem 1.3) because an affirmative answer would imply that the kernel of the Cassels' pairing equals the maximal divisible subgroup of  the Tate-Shafarevich group. This question appears again in \cite{cassels:IV} Problem (b) where Cassels completes his analysis of the kernel of Cassels' pairing but states that the question of the divisibility of $\Sha(A/k)$ in $\rH^1(k,A)$ remains open.

\smallskip

We will attempt to answer the above question one prime at a time. For a prime $p$ we say that $\Sha(A/k)$ is $p$-divisible in $\rH^1(k,A)$ if
\[
\Sha(A/k) \subseteq p^n \rH^1(k,A) \indent \text{for every }n\in \bN.
\] 
In view of the conjectured finiteness of $\Sha(A/k)$ which is known for elliptic curves over $\bQ$ of analytic rank $\leq 1$, the $p$-divisibility should be guaranteed for large $p$ depending on $A/k$.  Our aim therefore is to identify conditions that a prime number $p$ must satisfy for the $p$-divisibility conclusion to hold, which are sufficient and as close as possible to being necessary. 

\smallskip

For elliptic curves over $\bQ$, Theorem~\ref{thm:rationalEll} and 
Corollary~\ref{cor:quadratictwistoverQ} yield the following. 
\begin{thmA}
Let $E/\bQ$ be an elliptic curve defined over the rationals. Then the following holds.
\begin{enumerate}
\item  $\Sha(E/\bQ)$ is $p$-divisible in $\rH^1(\bQ,E)$ for all primes $p>7$.
\item There is at most one odd prime number $p = 3,5$ or $7$, and then at most two quadratic twists $E^\tau/\bQ$  of $E/\bQ$, such that $\Sha(E^\tau/\bQ)$ is not $p$-divisible in $\rH^1(\bQ,E^\tau)$. 
\item If $p=5$ or $7$ and $\Sha(E/\bQ)$ is not $p$-divisible in $\rH^1(\bQ,E)$ then $E$ has semistable reduction outside $p$.
\end{enumerate}
\end{thmA}

Note that our method can not handle the prime $p=2$ and Brendan Creutz \cite{Creutz} has now shown that for the elliptic curve 
$E : Y^2 = X(X + 80)(X + 205)$, labelled ``1025b2'' in Cremona's tables, $\Sha(E/\bQ)$ is not divisible by $4$ in $\rH^1(\bQ,E)$. More recently, Creutz found that divisibility by $9$ fails for the elliptic curve $E: X^3 + Y^3 + 138Z^3 = 0$. The remaining cases $p=5$ and $p=7$ for elliptic curves over $\bQ$ are open.

\smallskip

The essential content of Theorem~\ref{thm:numberfieldEll} provides a uniform criterion for elliptic curves over an algebraic number field $k$ that depends only on the degree of $k/\bQ$.  This is part (1) of Theorem~B below, see Section~\S\ref{sec:piofd} for the definition of $\pi(d)$; part (2) is proved as Corollary~\ref{cor:irreduciblecase}.

\begin{thmB}
Let $E/k$ be an elliptic curve defined over an algebraic number field $k$ of degree $d$ over $\bQ$. Then  
$\Sha(E/k)$ is $p$-divisible in $\rH^1(k,E)$ in the following two cases:
\begin{enumerate}
\item For every prime number $p> \max\{(2^d  + 2^{d/2})^2,\pi(d)\}$.
\item If $p \geq 3$, the $p$-torsion subgroup $E_p \subseteq E$ is an irreducible $\Gal_k$-representation, and $[k(\zeta_p):k] \not=2$ where $\zeta_p$ is a primitive $p$-th root of unity. 
\end{enumerate}
\end{thmB}

Using the bound  $\pi(d) < (1+3^{d/2})^2$ due to Oesterl\'e (unpublished), the bound in Theorem~B~(1) is simply $p>(2^d  + 2^{d/2})^2$. For more on $\pi(d)$ we refer to Section~\S\ref{sec:piofd}. 

\smallskip

Our method actually shows more. We define the \textbf{locally divisible $\rH^1$} as the kernel
\[
\rH^1_{\divisor}(k,A) = \ker\Big(\rH^1(k,A) \to \bigoplus_v \rH^1(k_v,A)/\divisor\big(\rH^1(k_v,A)\big) \Big)
\]
containing the global cohomology classes which locally become divisible. 
Since the formation of $\divisor(-)$ is a functor, we find 
\[
\divisor(\rH^1(k,A)) \subseteq \rH^1_{\divisor}(k,A).
\]
We now focus on the question of $p$-divisibility. Denote by $M_{p^n}$  the $p^n$-torsion of an abelian group $M$ and set $M_{p^\infty}= \cup_{n\in \bN} M_{p^n}$. By local Tate duality we have that
\[
\rH^1(k_v,A) = \Hom(A^t(k_v),\bQ/\bZ),
\]
where $A^t$ denotes the dual abelian variety of $A$. Then Mattuck's theorem on the structure of $A^t(k_v)$
implies that $\divisor(\rH^1(k_v,A)_{p^\infty}) = 0$ for $v \nmid p$. Hence, we have the following exact sequence
\begin{equation} \label{eq:pparth1localdiv}
0 \to \Sha(A/k)_{p^\infty} \to \rH^1_{\divisor}(k,A)_{p^\infty} \to \bigoplus_{v \mid p} \divisor\big(\rH^1(k_v,A)_{p^\infty}\big) \cong (\bQ_p/\bZ_p)^{\dim(A) \cdot [k:\bQ]}.
\end{equation}

The most general divisibility property that we prove can be viewed as a \textbf{local--to--global principle for $p$-divisibility} of elements in the Weil--Ch\^atelet group with respect to certain $p$. 

\begin{thmC}
Let $A/k$ be an abelian variety over an algebraic number field, $A^t$ its dual abelian variety, and $p$ a prime number. If we assume that 
\begin{enumerate}
\item[(i)] $\rH^1(k(A^t_p)/k,A^t_p) = 0$, with the splitting field $k(A^t_p)$ of the $p$-torsion $A^t_p$ of $A^t$, and 
\item[(ii)] the $\Gal_k$-modules $A^t_p$ and $\End(A^t_p)$ have no common irreducible subquotient,
\end{enumerate}
in particular if $p \gg 0$, then $\divisor(\rH^1(k,A))_{p^\infty} = \rH^1_{\divisor}(k,A)_{p^\infty}$  holds and $\Sha(A/k)$ is $p$-divisible in $\rH^1(k,A)$. 
\end{thmC}

In fact, Theorem A and Theorem B are, with the exception of two cases of elliptic curves $E/\bQ$ with respect to $p=11$, special cases of Theorem C. 
The proof of Theorem C and thus the answer to Cassels' divisibility question for almost all primes, actually follows immediately by Proposition \ref{prop:prdiv} from Theorem D  below (applied to 
$A^t$) and Theorem~\ref{thm:abelianalmostall}. 

\smallskip

Theorem D is actually a local global principle, i.e., a vanishing result for $\Sha^1(k,A_{p^n})$, see 
Section~\S\ref{sec:tassel} for the definition, that might be of independent interest; for the proof see Theorem~\ref{thm:proofofthmD}.

\begin{thmD} 
Let $A/k$ be an abelian variety over an algebraic number field, and let $p$ be  a prime number. 
Then 
\[
\Sha^1(k,A_{p^n}) = 0 \indent \text{ for all } n\geq 0,
\]
if we assume that 
\begin{enumerate}
\item[(i)] $\rH^1(k(A_p)/k,A_p) = 0$, with the splitting field $k(A_p)$ of the $p$-torsion $A_p$ of $A$, and 
\item[(ii)] the $\Gal_k$-modules $A_p$ and $\End(A_p)$ have no common irreducible subquotient.
\end{enumerate}
\end{thmD}

\smallskip

Our approach works also in the case when $A/k$ is an abelian variety over the function field of a geometrically connected, smooth projective curve $X$ over a finite field, at least with respect to divisibility questions for $p$ distinct from the residue characteristic. However, since our current interest lies in the arithmetic case when $k$ is an algebraic number field, we have not developed the geometric case in parallel.

%-----------------------------------------------------------------------------------------------------------------------

\smallskip

\subsection*{Acknowledgments} 
The authors would like to thank Brian Conrad, Brendan Creutz, Wojciech Gajda, and Joseph Oesterl\'e  for several useful discussions. The first author is also grateful to her advisor, Andrew Wiles, for introducing her to this method of thinking about the $p$-divisibility of the Tate-Shafarevich group. 

\smallskip

%-----------------------------------------------------------------------------------------------------------------------

\subsection{Notation}
We fix some notation which will be in use throughout the text.

\subsubsection{} Let $k$ denote an algebraic number field with absolute Galois group $\Gal_k$ and ring of integers $\fo_k$. 
The completion of $k$ at a place $v$ is $k_v$. The finite places of $k$ will be identified with the closed points of $X=\Spec(\fo_k)$. 

\subsubsection{} Let $M$ be an abelian group. The $n$-torsion subgroup is denoted by $M_n$, and $M_{p^\infty}$ denotes the $p$-primary torsion $\bigcup_n M_{p^n}$. The subgroup 
\[
\divisor(M)=\bigcap_{n\geq 1} nM
\]
of divisible elements of $M$ contains the maximal divisible subgroup of $M$
\[
\Div(M) = \sum_{\varphi \in \Hom ( \bQ , M)} \im ( \varphi ).
\]
For matters of clarity we stress that 
\textbf{"$a$ is $p$-divisible"} means that  for every $n \geq 1$ there is $a'_n$ with $a = p^n a'_n$.
Observe that $\divisor(M)$ can be strictly larger than $\Div(M)$, as for example with $d \in \bN$ in:
\[
M= \big(\bigoplus_n \frac{1}{dn}\bZ/\bZ\big)/\ker \big(\text{sum}:\bigoplus_n \frac{1}{d}\bZ/\bZ \rightarrow \frac{1}{d}\bZ/\bZ  \big)
\]
where $\divisor(M)=\frac{1}{d}\bZ/\bZ $ and $\Div(M)=0$.

\subsubsection{} \label{sec:piofd}
For every positive $d \in \bN$, we define $\pi(d)$ to be the maximal prime number $p$ such that there exists an elliptic curve $E$ defined over a number field $k$ of degree $d=[k:\bQ]$, which admits a non-trivial $k$-rational $p$-torsion point.

Merel \cite{merel:torsion} showed that $\pi(d)$ is finite and gave the upper bound ($d>1$)
\begin{equation} \label{eq:merelbound}
\pi(d) < d^{3d^2}.
\end{equation}
The best result in the literature to our knowledge  is the bound by Parent \cite{parent:bound}
\begin{equation} \label{eq:parentbound}
\pi(d) < 65\cdot (3^d-1)\cdot (2d)^6,
\end{equation}
while an unpublished result of Oesterl\'e\footnote{The improved bound for $\pi(d)$ due to Oesterl\'e from 1995 is unpublished. Moreover, Parent found a gap for $d=3$ and $p=43$, but was later able to repair the gap in  \cite{parent:torsion}  under an arithmetic assumption which was a consequence of the Birch and Swinnerton-Dyer conjecture. Fortunately, this assumption was later proved by Kato. We thank J.~Oesterl\'e for first hand information on his bound.} 
yields the estimate
\begin{equation} \label{eq:oesterlebound}
\pi(d) < (1+3^{d/2})^2.
\end{equation}
For small values of $d$ the following values of $\pi(d)$ are known precisely: Mazur in \cite{mazur:ratisog} Theorem~2
shows $\pi(1) = 7$, Kamienny \cite{kamienny:quadratic2} building on work of Kenku and Momose proves $\pi(2) = 13$, and Parent shows $\pi(3) = 13$ in \cite{parent:torsion} with the prime $17$ dealt with in \cite{parent:17}. Unpublished work by Kamienny, Stein and Stoll  (resp.\ together with Derickx) shows $\pi(4) = 17$ (resp.\ $\pi(5) = 19$).

%%%%%%%%%%%%%%%%%%%%%%%%%%%%%%%%%%%%%%%%%%

\section{Reminder on global Galois cohomology} 

%-----------------------------------------------------------------------------------------------------------------------

\subsection{Tate--Shafarevich groups and Poitou--Tate duality} \label{sec:tassel}

Let $k$ be an algebraic number field and $M$ a discrete $\Gal_k$-module. We set 
\[
\Sha^i(k,M) = \ker\Big(\rH^i(k,M) \xrightarrow{\prod_v \res_v} \prod_v \rH^i(k_v,M)\Big)
\]
with the restriction maps $\res_v$ induced by the embedding $k \inj k_v$, and the product ranges over all places $v$ of $k$.  

\smallskip

Let $M$ be a finite $\Gal_k$-module and $M^D  := \Hom(M,\mu_{\infty})$ its Cartier dual, here $\mu_{\infty}$ denotes the group of roots of unity in $k^\alg$. Poitou--Tate duality yields a perfect pairing of finite groups
\[
\Sha^1(k,M) \times \Sha^2(k,M^D) \to \bQ/\bZ,
\]
see \cite{nsw} VIII Theorem (8.6.7).

\smallskip

The Tate--Shafarevich group $\Sha(A/k)$ for an abelian variety $A$ over $k$ is defined\footnote{Indeed, the traditional definition of $\Sha(A/k)$ is equivalent due to the subtle equality
\[
\rH^1(k_v,A) = \rH^1(k_v,A(k_v^\alg)) = \rH^1(k_v,A(k^\alg)).
\]
} as 
\[
\Sha(A/k) = \Sha^1(k,A(k^\alg)),
\]
and in particular is a torsion group.  It follows that 
Cassels' question decomposes into $p$-primary parts.  We  will now concentrate on the $p$-primary part for a fixed prime number $p$.

%-----------------------------------------------------------------------------------------------------------------------

\subsection{The Selmer group and generalized Selmer groups}
Let $A/k$ be an abelian variety. 
The $p^n$-torsion Selmer group of $A$ is defined as
\[
\rH^1_\Sel(k,A_{p^n}) = \ker\Big(\rH^1(k,A_{p^n}) \to \prod_v \rH^1(k_v,A)\Big)
\]
with $v$ ranging over all places of $k$. A diagram chase with the cohomology sequence of the Kummer sequence $0 \to A_{p^n} \to A \xrightarrow{p^n\cdot} A \to 0$  over $k$ and all the localisations $k_v$ yields the fundamental short exact sequence
\begin{equation} \label{eq:selmerdef}
0 \to A(k)/p^nA(k) \xrightarrow{\delta_\kum} \rH^1_\Sel(k,A_{p^n}) \to \Sha(A/k)_{p^n} \to 0.
\end{equation}
It is known that $ \rH^1_\Sel(k,A_{p^n})$ is a finite group. 

\smallskip

The Selmer group $\rH^1_\Sel(k,A_{p^n})$ is a global $\rH^1$ with local Selmer-conditions determined by 
the image of the boundary map $\delta_{\kum}$ of the Kummer sequence 
\[
\Sel_v = \delta_\kum\big(A(k_v)/p^nA(k_v)\big) \subset \rH^1(k_v,A_{p^n})
\]
at every place. 
We shall be working with the following generalized Selmer groups. Let $Q$ be a finite set of finite places. The \textbf{Selmer group free at $Q$} is defined as 
\[
\rH^1_{\Sel_Q}(k,A_{p^n}) = \ker\Big(\rH^1(k,A_{p^n}) \to \prod_{v \not\in Q} \rH^1(k_v,A_{p^n})/\Sel_v\Big)
\]
and  the \textbf{Selmer group trivial at $Q$} is defined as 
\[
\rH^1_{\Sel^Q}(k,A_{p^n}) = \ker\Big(\rH^1(k,A_{p^n}) \to \prod_{v \not\in Q} \rH^1(k_v,A_{p^n})/\Sel_v \times \prod_{v \in Q} \rH^1(k_v,A_{p^n}) \Big).
\]
Consequently, we have the following inclusions 
\[
\rH^1_{\Sel^Q}(k,A_{p^n}) \subseteq \rH^1_{\Sel}(k,A_{p^n}) \subseteq \rH^1_{\Sel_Q}(k,A_{p^n}).
\]
In the case when $Q$ is the set of primes of $k$ dividing $p$, we use $\rH^1_{\Sel^p}(k,A_{p^n})$ (resp. $\rH^1_{\Sel_p}(k,A_{p^n})$) to denote the Selmer groups $\rH^1_{\Sel^Q}(k,A_{p^n})$ (resp. $\rH^1_{\Sel_Q}(k,A_{p^n})$).

%%%%%%%%%%%%%%%%%%%%%%%%%%%%%%%%%%%%%%%%%%
\section{Subgroups of \texorpdfstring{$\GL_2$}{GL-zwei} over finite prime fields}  \label{sec:group}
The purpose of this section is to provide the following classification statement.

\begin{thm} \label{thm:groupcrit}
Let $V$ be a vector space over $\bF_p$ of dimension $2$. For a subgroup $G \subseteq \GL(V)$ the following are equivalent.
\begin{enumerate}
\item 
\begin{enumerate}
\item[(i)] $V$ and $\End(V)$ have no common irreducible factor as $G$-modules, and
\item[(ii)] $\rH^1(G,V) = 0$.
\end{enumerate}
\item 
\begin{enumerate}
\item[(i)] The group $G$ is not contained in a subgroup of $\GL(V)$ that is isomorphic to the symmetric group $S_3$, in particular $p >2$, and 
 \item[(ii)] If $V$ is a reducible $G$-module, namely an extension 
 \begin{equation} \label{eq:extensionV}
 0 \to \chi_1 \to V \to \chi_2 \to 0
 \end{equation}
 for characters $\chi_i : G \to \bF_p^\ast$, then we require that $\chi_1 \not= \one , \chi_2^2$ and $\chi_2 \not= \one, \chi_1^2$.
\end{enumerate}
\end{enumerate}
\end{thm}

The proof of Theorem~\ref{thm:groupcrit} requires some preparation. 
We first recall the well known classification of subgroups of $\GL_2(\bF_p)$, see \S2 of \cite{Serre}.

\begin{prop} \label{prop:subgroups}
Let $G$ be a subgroup of $\GL_2(\bF_p)$ and let $\ov{G}$ be the image under the natural map $\GL_2(\bF_p) \to \PGL_2(\bF_p)$. Then one of the following  holds.
\begin{enumerate}
\item $p \mid \#G$ and $G$ is contained in a Borel $B \subset \GL_2(\bF_p)$.
\item $p \mid \#G$ and $G$ contains $\SL_2(\bF_p)$.
\item $p \nmid \#G$ and $G$ is contained in a normalizer of a split torus.
\item $p \nmid \#G$ and $G$ is contained in a normalizer of a non-split torus.
\item $p \nmid \#G$ and $\ov{G}$ is isomorphic to $A_4$, $S_4$, or $A_5$. \hfill $\square$
\end{enumerate}
\end{prop}

%----------------------------------------------------------------------------------------------------------------

\subsection{Computation of cohomology}

A $p$-Sylow subgroup $P$ of $\GL(V)$ is mapped under a suitable isomorphism $\GL(V) \cong \GL_2(\bF_p)$ induced by a choice of basis for $V$ to the standard $p$-Sylow subgroup
\[
U = \left\{\matzz{1}{b}{0}{1} \ ; \ b \in \bF_p\right\} \subset \GL_2(\bF_p).
\]
The normalizer $N$ of $P$ in $\GL(V)$ is mapped to the standard Borel subgroup
\[
B = \left\{\matzz{a_1}{b}{0}{a_2} \  ; \ a_1,a_2 \in \bF_p^\ast, \ b \in \bF_p\right\} \subseteq \GL_2(\bF_p).
\]
The two characters $\chi_i : N \to \bF_p^\ast$ defined by the isomorphism $N \cong B$ and 
\[
\chi_i(\matzz{a_1}{b}{0}{a_2}) = a_i,
\]
for $i=1,2$, occur in an exact sequence of $N$-modules
\begin{equation} \label{eq:asborelmodule}
0 \to \chi_1 \to V \to \chi_2 \to 0.
\end{equation}
The diagonal torus $T \subset B$ corresponds to a subgroup $S \subset N$ such that 
$N = P \rtimes S$ with $S$ acting on $P$ through the character 
\[
\chi_1 / \chi_2  : S \subset N \to \bF_p^\ast = \Aut(P).
\]

For a subgroup $G \subset \GL(V)$ with $p \mid \#G$, up to conjugation, we may assume that $P \subset G$ and determine the normalizer $N_G(P)$ of $P$ in $G$ as $N \cap G$. The extension (\ref{eq:asborelmodule}) restricts to an exact sequence as $N_G(P)$-modules in which by abuse of notation we write $\chi_i = \chi_i|_{N_G(P)}$. As before, the semi-direct product structure $N_G(P) = P \rtimes (S\cap G)$ has $S \cap G \cong N_G(P)/P$ acting on $P$ via the character $\chi_1/\chi_2$.

\begin{lem} \label{lem:computeH1}
(1) If $p \nmid \#G$, then $\rH^1(G,V) = 0$.

(2) If $p \mid \#G$, then the $N_G(P)$-equivariant quotient map  $\ph: V \surj \chi_2$ induces an injection
\[
\rH^1(G,V) \inj \Hom_{N_G(P)}(P,\chi_2).
\]
In particular, then $\rH^1(G,V) = 0$ except possibly if $\chi_1 = \chi_2^2$.
\end{lem}
\begin{pro}
(1) is clear and only recalled for completeness. For (2), as the index of $N_G(P)$ in $G$ is prime to $p$ we have $\rH^1(G,V) \inj \rH^1(N_G(P),V)$ so that we may assume $G = N_G(P)$. Since the index of $P$ in $G$ is prime to $p$, we find 
\[
\rH^1(G,V) = \rH^1(P,V)^{G/P}
\]
and it remains to show that the map $\rH^1(P,V) \to \rH^1(P,\chi_2) = \Hom(P,\chi_2)$ is injective. 
We consider the long exact cohomology sequence for the extension (\ref{eq:asborelmodule}). Since $\rH^0(P,\chi_1) = \rH^0(P,V)$, the connecting map
\[
\delta : \bF_p = \rH^0(P,\chi_2) \to \rH^1(P,\chi_1) = \Hom(P,\chi_1) \cong \bF_p
\]
is an isomorphism, and the map $\rH^1(P,V) \to \rH^1(P,\chi_2)$ is indeed injective.
\end{pro}

%----------------------------------------------------------------------------------------------------------------

\subsection{When there are homotheties}

The center of $\GL(V)$ is the group $Z \cong \bF_p^\ast$ of scalar automorphisms. We formulate a more general lemma, because later our proof of Theorem~\ref{thm:abelianalmostall} depends on it.

\begin{lem} \label{lem:center} Let $W$ be a finite dimensional $\bF_p$-vector space, and let $G \subset \GL(W)$ be a subgroup which intersects the center $\bF_p^\ast \cong Z \subset \GL(W)$ non-trivially. Then the following holds.
\begin{enumerate}
\item $W$ and the adjoint representation $\End(W)$ have no common irreducible factor.
\item $\rH^1(G,W) = 0$.
\end{enumerate}
\end{lem}
\begin{pro}
(1) The group $H=G \cap Z$ of homotheties in $G$ acts trivially on every irreducible factor of $\End(W)$ and faithfully on every irreducible factor of $W$. Hence none of them can occur in both $W$ and $\End(W)$.

(2) The inflation/restriction sequence for $H \lhd G$ reads 
\[
0 \to \rH^1(G/H,W^H) \to \rH^1(G,W) \to \rH^1(H,W)^{G/H}.
\]
Since $H$ was assumed to be non-trivial and is necessarily of order prime to $p$, both $W^H$ and $ \rH^1(H,W)$ vanish, and consequently also $ \rH^1(G,W) = 0$. 
\end{pro}

\begin{lem} \label{lem:exceptionalcase}
Let $G$ be a subgroup of $\GL_2(\bF_p)$ such that (5) of Proposition~\ref{prop:subgroups} holds. Then $G$ meets the center $Z$ of $\GL_2(\bF_p)$ non-trivially.
\end{lem}
\begin{pro}
It suffices to discuss the case $\ov{G} = A_4$. If $G \cap Z = 1$, then we have a copy 
\[
A_4 \subseteq \GL_2(\bF_p),
\]
in particular $p >2$. The $2$-Sylow subgroup of $A_4$, the Klein $4$-group $V_4 = \bZ/2\bZ \times \bZ/2\bZ$, has a completely reducible representation theory already rationally over $\bF_p$ as we may produce enough projectors already rationally over $\bF_p$. Hence $V_4$ is contained in a split torus $C = \bF_p^\ast \times \bF_p^\ast$ and must agree with the $2$-torsion of $C$.   Thus $V_4$ already contains the central element 
\[
-1 \in \bF_p^\ast \cong Z,
\]
a contradiction.
\end{pro}

%----------------------------------------------------------------------------------------------------------------

\subsection{The adjoint representation}

The action of the normalizer of a torus is best understood by identifying $V = A = \bF_p[\alpha]$ with a separable quadratic $\bF_p$-algebra. Then 
\[
\Aut(A/\bF_p) = \bZ/2\bZ,
\]
of which we denote the generator by $F$ (since it is the Frobenius if $A$ is a field).  Let us consider  the associated  normalizer of a torus $G = A^\ast \rtimes \langle F \rangle \subset \GL(V)$. We introduce two new $G$-module structures on $V$ denoted 
\begin{itemize}
\item $A_0$
with action factoring through $G \surj \langle F \rangle$ followed by the $F$-action on $V=A_0$, and
\item $A_1$
with action factoring through
\[
G= A^\ast \rtimes \langle F \rangle \to    A^\ast \rtimes \langle F \rangle  = G
\]
which maps  $F$ to $F$ and $\lambda \in A^\ast$ to $\lambda/F(\lambda)$,
followed by the defining action on $V = A_1$.
\end{itemize}

\begin{lem} \label{lem:adjointfornormalizertorus}
In the above notation, $\End(V)$ decomposes as a $G = \bF_p[\alpha]^\ast \rtimes \langle F \rangle$-module as 
\[
A_0  \oplus A_1 \xrightarrow{\sim} \End(V)
\]
\[
(a_0,a_1) \mapsto \big(v \mapsto a_0 v + a_1 F(v)\big).
\]
\end{lem}
\begin{pro}
We compute compatibility with $F$
\[
(F(a_0),F(a_1)) \mapsto  \big(v \mapsto F(a_0) v + F(a_1) F(v)\big) = F \circ  \big(v \mapsto a_0 v + a_1 F(v)\big) \circ F^{-1}
\]
and compatibility wtih $\lambda \in A^\ast$
\[
(a_0, \frac{\lambda}{F(\lambda)} \cdot a_1) \mapsto  \big(v \mapsto a_0 v +  \frac{\lambda}{F(\lambda)} a_1 F(v)\big) =  \big(v \mapsto \lambda \cdot (a_0 (\lambda^{-1} v) + a_1 F(\lambda^{-1}\cdot v))\big).
\]

Clearly, the map restricted to either $A_i$ is injective. Since $A^\ast$ acts trivially on $A_0$ but non-trivially on $A_1$ (except for the case $A = \bF_2 \times \bF_2$ were a direct computation confirms the lemma) this implies that the sum of the map still is injective. By comparing dimensions we deduce that it is an isomorphism.
\end{pro}

\begin{lem} \label{lem:adjointforborel}
Let $V$ be a reducible representation with characters $\chi_i:G \to \bF_p^\ast$ for $i=1,2$ as irreducible factors. Then the irreducible factors of $\End(V)$ are $\chi_1 \otimes \chi_2^{-1}, \one, \one$, and $\chi_2 \otimes \chi_1^{-1}$.
\end{lem}
\begin{pro}
Choosing the right ordering we have a short exact sequence $0 \to \chi_1 \to V \to \chi_2 \to 0$. The dual $V^\vee$ then sits in the short exact sequence $0 \to \chi_2^{-1} \to V^\vee \to \chi_1^{-1} \to 0$, and tensoring the two sequences yields the desired decomposition of $\End(V) = V^\vee \otimes V$.
\end{pro}

\medskip

We next address the special case where $G \cong S_3$.

\begin{lem} \label{lem:S3}
If $G \cong S_3$, then $V$ and $\End(V)$ have a common irreducible factor. 
\end{lem}
\begin{pro}
We discuss the cases $p=2$, $p=3$ and $p \nmid \#G$ separately.

\smallskip

\textit{The case $p=2$.} Here a good model is $V = \bF_4$ and $G = \GL_2(\bF_2)$ acts as $\bF_4^\ast \rtimes \Gal(\bF_4/\bF_2)$. As in Lemma~\ref{lem:adjointfornormalizertorus} we denote the generator of $\Gal(\bF_4/\bF_2)$ in $\GL(V)$ by $F$. By Lemma~\ref{lem:adjointfornormalizertorus} the map $V \to A_1$ which maps $v \mapsto F(v)$ is an isomorphism of $V$ onto a direct summand of $\End(V)$. Indeed, for $\lambda \in \bF_4^\ast$ we find
\[
F(\lambda v) = \lambda^2 F(v) = \lambda^{-1} F(v) 
= \bruch{\lambda}{F(\lambda)} F(v)
\]
while the compatibility with $F$ is obvious.

\smallskip

\textit{The case $p=3$.} Here the $3$-cycle of $G=S_3$ must act after choosing a suitable basis as 
\[
\matzz{1}{1}{0}{1}.
\]
Therefore $G$ is contained in the standard Borel. As $S_3$ is not abelian, the two characters $\chi:G \to \bF_3^\ast$ by projecting to a diagonal entry must be distinct. Hence one is trivial and the other one is the sign character.  By Lemma~\ref{lem:adjointforborel} the trivial character also occurs in $\End(V)$.

\smallskip

\textit{The case $p \nmid \#G$.}  
By representation theory prime to $p$ there is a unique faithful representation of dimension $2$ of $G=S_3$ over $\bF_p$. Namely, the $3$-cycle lies in a torus $C$, either as 
\[
\matzz{\zeta_3}{0}{0}{\zeta_3^{-1}}
\]
in the split torus if $3 \mid p-1$ (here $\zeta_3$ denotes a primitive cubic root of unity), or as the $3$-torsion of a nonsplit torus if $3 \mid p+1$. With the notation of Lemma~\ref{lem:adjointfornormalizertorus}, a transposition agrees with the generator $F$ of a spliting $N/C \to N$ for the normalizer $N$ of the torus.  The action on $A_1 \subset \End(V)$ is also faithful by the formula  from Lemma~\ref{lem:adjointfornormalizertorus} as otherwise $F$ would fix the $3$-cycle and $G=S_3$ would be commutative, a contradiction. The uniqueness of a faithful $2$-dimensional representation for $S_3$ thus shows that $V \cong A_1 \subset \End(V)$ as $G$-representations.
\end{pro}

\begin{cor} \label{cor:detS3}
The determinant of the unique faithful $2$-dimensional $S_3$ representation $V$ over $\bF_p$ is the mod $p$ sign character ${\rm sign} : S_3 \to \{\pm 1\} \subseteq \bF_p^\ast$.
\end{cor}	
\begin{pro}
This follows from the discussion of this representation in the proof or Lemma~\ref{lem:S3}. Alternatively, we can remark that the representation V is already defined over $\bZ$ as the kernel $V_{\bZ}$ of the natural map 
\[
{\rm sum} \ : \ \ind_{(2,3)}^{S_3} (\one) \to \one
\]
which sums up the components in the natural basis of the tautological permutation representation. Then \[
\det(V_{\bZ}) = \det(\ind_{(2,3)}^{S_3} (\one)) = {\rm sign},
\]
according to one of the definitions of the sign character.
\end{pro}

%----------------------------------------------------------------------------------------------------------------

\subsection{Proof of Theorem~\ref{thm:groupcrit}}
By Lemma~\ref{lem:S3}, both properties (1) and (2) fail if $G$ is contained in a subgroup of $\GL(V)$ that is isomorphic to $S_3$. We thus may assume (i) of (2), and in particular that $p \geq 3$, since $\GL_2(\bF_2) \cong S_3$.

\smallskip

\textit{The case $V$ reducible.} Let $V$ as a $G$-module be an extension of the character $\chi_2$ by the character $\chi_1$.  By Lemma~\ref{lem:adjointforborel}, the irreducible factors of $\End(V)$ as a $G$-module are $\chi_1 \otimes \chi_2^{-1}, \one,$ and $\chi_2 \otimes \chi_1^{-1}$. So (i) of (1) holds if and only if $\chi_1 \not= \one, \chi_2^2$ and $\chi_2 \not= \one, \chi_1^2$, which is exactly what (ii) of (2) asks for. Moreover, Lemma~\ref{lem:computeH1} shows that (ii) of (1) follows from (ii) of (2).

\smallskip

\textit{The case $V$ irreducible.}
Let $p \mid \#G$. As $G$ is not contained in a Borel, we conclude by Proposition~\ref{prop:subgroups} that $G$ contains $\SL_2(\bF_p)$. Now as $p \geq 3$ the group $G$ necessarily meets the center of $\GL(V)$ non-trivially, so that (1) holds by Lemma~\ref{lem:center}. 

If $p \nmid \#G$ and the image of $G$ in $\PGL_2(\bF_p)$ is one of the exceptional cases $A_4$, $S_4$ or $A_5$, namely in case (5) of Proposition~\ref{prop:subgroups}, then Lemma 
\ref{lem:exceptionalcase} and Lemma~\ref{lem:center} show that (1) holds.  

It remains to show (1) in the case that $p \nmid \#G$ and $G$ is contained in the normalizer $N = C \rtimes \bZ/2\bZ$ of a torus $C \subset \GL_2(\bF_p)$ and $V$ is an irreducible $G$-module. Part (ii) of (1) holds trivially. Lemma~\ref{lem:adjointfornormalizertorus} now tells us about the action of $G$ on $\End(V)$. Let $H = G \cap C$ be the intersection with the torus, hence a subgroup in $G$ of index $\leq 2$. Because $V$ is not a reducible $G$-module and $p \geq 3$ we have $\# G \geq 3$ and therefore $H \not=1$. 

If (1) part (i) fails, then in the decomposition $\End(V) = A_0 \oplus A_1$ of Lemma~\ref{lem:adjointfornormalizertorus}  we must have $V \cong A_1$, because $H$ acts trivially on the first summand $A_0 \subset \End(V)$. The representation $V \otimes \bF_{p^2}$ regarded as an $H$-module decomposes as a
sum of characters $\chi_1 \oplus \chi_2$ of $H$. The representation $A_1 \subset \End(V)$ decomposes after scalar extension to $\bF_{p^2}$ as $H$-module as $\chi_1 \chi_2^{-1} \oplus \chi_2 \chi_1^{-1}$. Comparing the two, we find either $\chi_1 = \one = \chi_2$, whence $H=1$ contradicting the irreduciblity of $V$ as a $G$-module. Or, $\chi_1  =\chi_2^2 \not=1$ and $\chi_2 = \chi_1^2\not=1$ which means $\chi_1$ and $\chi_2$ are of order $3$ and determine each other. In this case $H \cong \bZ/3\bZ$  and acts on $V$ non-centrally and without a common fixed vector. In any case, split or non-split, the subgroup $H \subset C$ is normal but not central in $N$. Hence  either $H = G$ and $G$ can be embedded in a subgroup $S_3$ of $\GL(V)$, or $H \lhd G$ of index $2$ and $G \cong S_3$ itself. In any case, this violates (i) of (2) and was excluded in the beginning of the proof. This completes the proof of Theorem~\ref{thm:groupcrit}.
\hfill $\square$

%----------------------------------------------------------------------------------------------------------------
\subsection{Families of exceptional cases in \texorpdfstring{$\GL_{n}$}{GLn}}
In Theorem~\ref{thm:groupcrit}, the family of subgroups
\[
S_3 \subseteq \GL_2(\bF_p)
\]
provides a family of exceptions, which according to the theorem is the only family (varying $p$) acting irreducibly in dimension $n=2$ such that condition (2)  fails. We would like to illustrate that this is only the tip of the iceberg when the dimension $n$ grows.

\begin{ex} \label{ex:autalgebra}
Let $G/\bQ$ be a not necessarily connected linear algebraic group. By \cite{gordeevpopov:autofalgebra} Theorem~1, there is a simple finite dimensional algebra $A/\bQ$, in general non-associative and non-commutative. Hence, we have a non-trivial bilinear map 
\[
m: A \otimes_{\bQ} A \to A,
\]
such that the $\bQ$-algebraic group $\underline{\Aut}(A)$ is isomorphic over $\bQ$ to $G$. We regard $A$ as a (faithful) algebraic $G$-representation 
\begin{equation} \label{eq:reprautalgebra}
G = \underline{\Aut}(A) \inj \GL(A).
\end{equation}
The adjoint to $m$ with respect to the adjoint pair of functors $-\otimes_\bQ A$ and $\Hom_{\bQ}(A,-)$ on $G$-representations
\[
m^{\#} : A \to \Hom_{\bQ}(A,A) = \End(A)
\]
is a non-trivial map of $G$-representations.

If we choose integral models for all structures involved, we can reduce mod $p$ for almost all $p$ and find another exceptional family of subgroups 
\[
G(\bF_p) \subseteq \GL_n(\bF_p)
\]
with $n = \dim_\bQ(A)$ such that Theorem~\ref{thm:groupcrit} (2) fails.
\end{ex}

It is not clear to us, whether the representations \eqref{eq:reprautalgebra} in Example~\ref{ex:autalgebra} are necessarily irreducible. However, also when we assume the representation to be irreducible, we are not spared of more exceptional cases.

\begin{ex}
The monster finite simple group $M$ is the largest sporadic simple group of order
\[
|M| = 2^{46} \cdot 3^{20} \cdot 5^9 \cdot 7^6 \cdot 11^2 \cdot 13^3 \cdot 17 \cdot 19 \cdot 23 \cdot 29 \cdot 31 \cdot 41 \cdot 47 \cdot 59 \cdot 71 \approx 8 \cdot 10^{53}.
\]
It was predicted by Fischer and independently Griess in 1973 and later constructed by Griess \cite{griess:monster} as an automorphism group of the (real) Griess algebra 
\[
A = \one \oplus V
\]
of  dimension $1 + 196883$. The action on the $196883$-dimensional piece $V$ yields the smallest irreducible representation of $M$. The $M$-action on $V$ preserves the non-trivial product 
\[
\mu : V \otimes V  \subset A \otimes A \to A \to V
\]
that is the $V$-component  of the product of the Griess algebra $A \otimes A \to A$ projected to $V$. This product is defined over a ring $R$ of finite type over $\bZ$ and thus admits non-trivial specialisations to fields $\bF_p$ for a positive Dirichlet density of primes $p$. In fact, by \cite{griess:monster} page 2, the irreducible character $\chi$ of $V$ takes rational values, and therefore by the Brauer--Speiser Theorem \cite{curtisreiner:2} (74.27) has Schur index equal to $1$. It can thus be defined over $\bQ$ by \cite{curtisreiner:2} Theorem (74.5) (iii), and consequently, a non-trivial multiplication 
\[
\mu : V \otimes V \to V
\]
can even be defined over $R=\bQ$. Thus all but finitely many primes are ok. This shows that for dimension $n=196883$ there is a series of exceptional subgroups 
\[
M \subset \GL_n(\bF_p)
\]
isomorphic to the monster group, where Theorem~\ref{thm:groupcrit} (2) fails.
\end{ex}

%%%%%%%%%%%%%%%%%%%%%%%%%%%%%%%%%%%%%%%%%%

\section{Vanishing of the Selmer group trivial at \texorpdfstring{$Q$}{Q}} \label{sec:induction}

We recall that $k$ will always be an algebraic number field with ring of integers $\fo_k$. 

%----------------------------------------------------------------------------------------------------------------
\subsection{The obstruction for a local--global principle for divisibility}
Divisibility by $p^n$ is controlled by an obstruction class as follows.

\begin{lem}  \label{lem:probstruction}
Let $A/k$ be an abelian variety  and let $p$ be a prime number.
For $n \in \bN$ there is a short exact sequence
\begin{equation} \label{eq:delta_r}
0 \to  \rH^1_{\divisor}(k,A)_{p^\infty} \cap p^n \cdot \rH^1(k, A) \to  \rH^1_{\divisor}(k,A)_{p^\infty} \xrightarrow{\delta_n} 
\Sha^2(k,A_{p^n}).
\end{equation}
\end{lem}

\begin{pro}
The short exact sequence $0 \to A_{p^n} \to A \xrightarrow{p^n \cdot} A \to 0$ of $\Gal_k$-modules together with restriction to $\Gal_{k_v}$ yields a commutative diagram 
\[
\xymatrix@M+1ex@R-2ex{
 \rH^1(k, A) \ar[r]^{p^n \cdot}  \ar[d] &  \rH^1(k, A) \ar[r]^{\delta_n} \ar[d] &   \rH^2(k, A_{p^n}) \ar[d] \\
  \rH^1(k_v, A) \ar[r]^{p^n \cdot} &  \rH^1(k_v, A) \ar[r]^{\delta_n} &   \rH^2(k_v, A_{p^n}) .}
\]
As classes in $\rH^1_{\divisor}(k,A)$ are locally divisible, a diagram chase 
shows that $\delta_n(\rH^1_{\divisor}(k,A))$ takes values in 
$\Sha^2(k,A_{p^n})$. This proves the lemma.
\end{pro}

\begin{prop} \label{prop:prdiv} 
Let $A/k$ be an abelian variety and let $p$ be a prime number. If for $n \in \bN$ with the dual abelian varity $A^t$ we have $\Sha^1(k,A^t_{p^n}) = 0$, then $\rH^1_{\divisor}(k,A)_{p^\infty} \subseteq p^n \cdot \rH^1(k,A)$.
\end{prop}
\begin{pro}
By Lemma~\ref{lem:probstruction} and Poitou--Tate duality, see \cite{nsw} VIII Theorem (8.6.7),  the obstruction for divisibility by $p^n$ lies in  
\[
\Sha^2(k,A_{p^n}) = \Hom(\Sha^1(k,A^t_{p^n}),\bQ/\bZ) = 0.
\]
Now the proposition follows immediately from the exact sequence \eqref{eq:delta_r} of Lemma~\ref{lem:probstruction}.
\end{pro}

\smallskip

Let $Q$ be a finite set of primes. Then obviously
\[
\Sha^1(k,A^t_{p^n}) \subseteq \rH^1_{\Sel^Q}(k,A^t_{p^n})
\]
so that in view of Proposition~\ref{prop:prdiv} we can deduce divisibility by $p^n$ from vanishing theorems for Selmer groups trivial at $Q$.

%----------------------------------------------------------------------------------------------------------------
\subsection{The obstruction for divisibility in Cassels' question}

As far as Cassels' original question is concerned, namely the divisibility by $p^n$ of  only $\Sha(A/k)$ instead of $\rH^1_{\divisor}(k,A)$ in $\rH^1(k,A)$, we can formulate a necessary and sufficient condition\footnote{The statement of Proposition \ref{prop:casselsiff} was suggested to us by Brendan Creutz.} as follows.

\begin{prop} \label{prop:casselsiff}
Let $A/k$ be an abelian variety and let $p$ be a prime number. For $n \in \bN$ the following are equivalent.
\begin{enumerate}
\item[(i)] $\Sha(A/k) \subseteq p^n \cdot \rH^1(k,A)$.
\item[(ii)] The natural map $i_\ast : \Sha^1(k,A^t_{p^n}) \to \Sha(A^t/k)$ factors over $\Div\big(\Sha(A^t/k)\big)$.
\end{enumerate}
\end{prop}
\begin{pro}
By Lemma~\ref{lem:probstruction} assertion (i) is equivalent to the vanishing of the Kummer boundary map restricted to $\Sha(A/k)$ which factors as 
\begin{equation} \label{eq:deltasha}
\delta : \Sha(A/k)/\Div(\Sha(A/k)) \to \Sha^2(k,A_{p^n}).
\end{equation}
Assertion (ii) is clearly equivalent to the vanishing of 
\begin{equation} \label{eq:inclusionsha}
\Sha^1(k,A^t_{p^n}) \to \Sha(A^t/k)/\Div(\Sha(A^t/k)),
\end{equation}
so that it remains to argue that \eqref{eq:inclusionsha} is the adjoint of \eqref{eq:deltasha} under the Poitou--Tate pairing $\langle-,-\rangle_{\rm PT}$ recalled in Section \S\ref{sec:tassel} and the Cassels--Tate pairing $\langle-,-\rangle_{\rm CT}$, since the latter is non-degenerate modulo the maximal divisible subgroup.

For $\alpha' \in \Sha^1(k,A^t_{p^n})$ and $\alpha \in \Sha(A/k)$ we compute the Cassels--Tate pairing of $i_\ast(\alpha')$ with $\alpha$ following basically the Weil-pairing definition of the Cassels--Tate pairing as in \cite{PS} \S12.2. Instead of computing the boundary of $\alpha$ after choosing $p^n$ large enough and lifting $\alpha$ to $\rH^1(k,A_{p^n})$ with respect to 
\[
0 \to A_{p^n} \to A_{p^{2n}} \xrightarrow{p^n} A_{p^n} \to 0,
\]
we can similarly compute the boundary directly with respect to 
\[
0 \to A_{p^n} \to A \xrightarrow{p^n} A \to 0.
\]
Therefore $\langle i_\ast(\alpha'), \alpha \rangle_{\rm CT}$ depends only on $\alpha'$ and $\delta(\alpha)$, and the further formulae in 
\cite{PS} \S12.2 translate immediately into the description of the Poitou--Tate pairing $\langle \alpha', \delta(\alpha) \rangle_{\rm PT}$
as is explicitly described in the paragraph above Theorem $4.20$ in \cite{Milne}. This finishes the proof. 
\end{pro}

\begin{rmk}
(1)
In this paper, assertion (ii) of Proposition~\ref{prop:casselsiff} will always follow from the vanishing of $\Sha^1(k,A_{p^n})$. Therefore, in all cases the more general assertion on divisibility properties of $\rH^1_{\divisor}(k,A)$ will follow.

(2) 
It follows from \eqref{eq:selmerdef} that we have a short exact sequence
\begin{equation} \label{eq:locdivpointssequence}
0 \to \{a \in A(k) ; a \in p^n \cdot A(k_v) \text{ for all $v$ }\}/p^n \cdot A(k) \to \Sha^1(k,A_{p^n}) \to \Sha(A/k)_{p^n}
\end{equation}
Cases where $\{a \in A(k) ; a \in p^n \cdot A(k_v) \text{ for all $v$ }\} \not= p^n \cdot A(k)$ are known, even when $k=\bQ$ and $A=E$ is an elliptic curve: see \cite{dz:4} for $p^n=4$, and  more generally in \cite{paladino:counterex}. These examples yield in particular cases where 
\[
\Sha^1(k,A_{p^n}) \not= 0,
\]
but in view of Proposition~\ref{prop:casselsiff} and \eqref{eq:locdivpointssequence} this does not imply that Cassels' question has a negative answer for $A^t/k$. 

(3) 
Brendan Creutz in \cite{Creutz} also proves Proposition~\ref{prop:casselsiff} and uses it to show that for every prime $p$ there exists an abelian variety $A/\bQ$ such that $\Sha(A/\bQ)$ is not divisible by $p$ in $ \rH^1(\bQ,A)$. Furthermore, he also finds an elliptic curve $E/\bQ$ such that $\Sha(E/\bQ)_2$ is not divisible by $4$ in $\rH^1(\bQ,E)$.
\end{rmk}

%----------------------------------------------------------------------------------------------------------------
\subsection{The Selmer splitting field}   \label{sec:Selmersplitting}
Let $A/k$ be an abelian variety.  We denote by
\[
\rH^1_{\Sel}(k(A_p)/k,A_p)
\]
the intersection of $\rH^1_{\Sel}(k, A_p)$ with the image under inflation of $\rH^1(k(A_p)/k,A_p)$  in $\rH^1(k,A_p)$.
Then we have the following commutative diagram with exact rows 
\[ 
\xymatrix@M+1ex@R-2ex{
0 \ar[r] &  \rH^1_{\Sel}(k(A_p)/k,A_p) \ar@{}[d]|{{\rotatebox{-90}{$\subseteq$}}} \ar[r]^(0.55){\infl} & \rH^1_{\Sel}(k, A_p) \ar[r]^(0.3){\res}  \ar@{}[d]|{{\rotatebox{-90}{$\subseteq$}}} & \Hom_{\Gal(k(A_p)/k)}(\Gal_{k(A_p)}^\ab \otimes \bF_p,A_p)  \ar@{=}[d] \\
0 \ar[r] &  \rH^1(k(A_p)/k,A_p) \ar[r]^(0.55){\infl} & \rH^1(k, A_p) \ar[r]^(0.3){\res} & \rH^1(k(A_p),A_p)^{\Gal(k(A_p)/k)} 
} 
\] 
The restriction map defines a canonical continuous $\Gal_k$-equivariant pairing 
\begin{equation} \label{pairing-selmer}
 \rH^1_{\Sel}(k, A_p) \times \big(\Gal_{k(A_p)}^\ab \otimes \bF_p\big) \to A_p.
\end{equation}
Let $M$ denote the quotient of $\Gal_{k(A_p)}^\ab \otimes \bF_p$ by the right kernel of the above pairing. The restriction map factors as follows
\[
\rH^1_{\Sel}(k, A_p) \to  \Hom_{\Gal(k(A_p)/k)}(M,A_p)  \subset   \Hom_{\Gal(k(A_p)/k)}(\Gal_{k(A_p)}^\ab \otimes \bF_p,A_p) .
\]
Observe that the finiteness of the Selmer group $ \rH^1_{\Sel}(k, A_p)$ implies that $M$ is finite. Hence the quotient $M$ corresponds to a finite Galois extension $L/k(A_p)$,  that we call the \textbf{Selmer splitting field} of $A$ with respect to the prime $p$, more precisely $M=\Gal(L/k(A_p))$. Since $M$ is a quotient as $\Gal(k(A_p)/k)$-module, the field $L$ is in fact Galois over $k$ and $\Gal(k(A_p)/k)$ acts on $M$ by conjugation after lifting under the quotient map $\Gal(L/k) \surj \Gal(k(A_p)/k)$. 

\begin{lem} \label{lem:selmerslpittingfield}
Let $A/k$ be an abelian variety, and let $L$ be the Selmer splitting field with respect to $p$. Then  the following holds.
\begin{enumerate}
\item The following sequence is exact:
\begin{equation} \label{eq:selmersplittingfield}
0 \to \rH^1_{\Sel}(k(A_p)/k,A_p) \to \rH^1_{\Sel}(k, A_p) \to  \Hom_{\Gal(k(A_p)/k)}(\Gal(L/k(A_p)),A_p)
\end{equation}
\item Every irreducible $\Gal(k(A_p)/k)$-module subquotient of $\Gal(L/k(A_p))$ is isomorphic to an irreducible subquotient of $A_p$.
\end{enumerate}
\end{lem}
\begin{pro}
(1) follows from the definition of $L$. For (2) we note that pairing \ref{pairing-selmer} yields an injective map
\begin{equation} \label{eq:structureM}
\Gal(L/k(A_p)) = M \inj \Hom( \rH^1_{\Sel}(k, A_p),A_p) \cong A_p \oplus \ldots \oplus A_p
\end{equation}
of  $\Gal(k(A_p)/k)$-modules, where $\rH^1_{\Sel}(k, A_p)$ carries the trivial action.
\end{pro}

%----------------------------------------------------------------------------------------------------------------
\subsection{The vanishing of some generalized Selmer groups}  

\begin{prop}  \label{prop:structureofgenselmer}
Let $A/k$ be an abelian variety, and let $p$ be a prime number, such that
\begin{enumerate}
\item[(i)] $\rH^1(k(A_p)/k,A_p) = 0$.
\end{enumerate}
Let $Q$ be a finite set of finite primes of $k$ not dividing $p$, and fix $n \in \bN$ such that 
\begin{enumerate}
\item[(ii)] $A$ has good reduction at $v$ for all $v \in Q$;
\item[(iii)]  the set of Frobenius elements $\Frob_{w} \in \Gal(L/k(A_p))$ where $L$ is the Selmer splitting field of $A/k$ with respect to $p$ and $w$ denotes a prime of $k(A_p)$ dividing $v$, when $v$ ranges over $Q$,  generates $\Gal(L/k(A_p))$ as a $\Gal(k(A_p)/k)$-module;
\item[(iv)] $A_{p^n}(k_v)$ is a free $\bZ/p^n\bZ$-module for all $v \in Q$.
\end{enumerate}
Then for all $m \leq n$ we have that
\[
\rH^1_{\Sel^Q}(k,A_{p^m}) = 0.
\]
\end{prop}

\begin{pro}
\textit{Step 1:}  We first treat  $m=1$.  We set $k_1 = k(A_p)$, and $k_{1,w}$ for the completion of $k(A_p)$ in $w$. Localization at $v$ yields a commutative diagram
\[
\xymatrix@M+1ex@R-1ex@C+2ex{ \rH^1_{\Sel}(k, A_p) \ar[r]^(0.35){\res_{k_1/k}} \ar[d] & \Hom_{\Gal(k_1/k)}(\Gal(L/k_1),A_p) \ar[d]^{\ev_w} \\
\rH^1_{\nr}(k_v,A_p) \ar[r]^(0.45){\res_{k_{1,w}/k_v}} & \rH^1_{\nr}(k_{1,w},A_p) = A_p}
\]
with the evaluation map $\ev_w$ mapping a morphism $\ph : \Gal(L/k_1)  \to A_p$ to its value $\ph(\Frob_{w})$ at the Frobenius element of $w$. Assumption (iii), the sequence \eqref{eq:selmersplittingfield}, and assumption (i) imply
\[
\rH^1_{\Sel^Q}(k,A_p) \subseteq \rH^1_{\Sel}(k(A_p)/k,A_p) \subseteq \rH^1(k(A_p)/k,A_p) = 0.
\]

\smallskip

\textit{Step 2:} We now show the general case by induction on $m$ terminating in $n$. As an abbreviation we set 
\[
\bL_{m,v} = \Sel^Q_v \subseteq \rH^1(k_v,A_{p^m})
\]
for the Selmer condition trivial at $Q$ for $A_{p^m}$-coefficients. Then 
 the following diagram is commutative and the rows are exact
\begin{equation} \label{eq:compatiblelocalconditions}
\xymatrix@M+1ex@R-2ex{   
 & {\bL_{m-1,v}} \ar[r] \ar[d] &  {\bL_{m,v}} \ar[r] \ar[d]  &  {\bL_{1,v}} \ar[d] \ar[r] & 0  \\
0 \ar[r] & \rH^1(k_v,A_{p^{m-1}})/\delta_\kum\big(A_p(k_v)\big) \ar[r] &  \rH^1(k_v,A_{p^m}) \ar[r]^{p^{m-1}\cdot} &  \rH^1(k_v,A_{p})
}
\end{equation}
The snake lemma applied to \eqref{eq:compatiblelocalconditions} yields that in the commutative diagram
\begin{equation} \label{eq:almostSel}
\xymatrix@M+1ex@R-2ex@C-2ex{ A_p(k) \ar[d] \ar[r]^(0.4){\delta_\kum} &  \rH^1(k,A_{p^{m-1}})  \ar[d] \ar[r] & \rH^1(k,A_{p^m}) \ar[r]  \ar[d] & \rH^1(k,A_p)  \ar[d] \\
   \prod_{v} A_p(k_v) \ar[r]^(0.4){\delta_\kum} &  \prod_{v } \bruch{\rH^1(k_v,A_{p^{m-1}})}{\bL_{m-1,v}} \ar[r] & \prod_{v} \bruch{\rH^1(k_v,A_{p^m})}{\bL_{m,v}} \ar[r] & \prod_{v} \bruch{\rH^1(k_v,A_p)}{\bL_{1,v}}
}
\end{equation} 
the bottom row is exact.  The map $\delta_\kum$ in the bottom row is the zero map: by assumption (iv) when $v \in Q$, then $A_{p^m}(k_v)/A_{p^{m-1}}(k_v) \to A_p(k_v)$ is surjective for $m \leq n$, or, in general for $v \notin Q$, by comparing the boundary maps for the diagram 
\[
\xymatrix@M+1ex@R-2ex{
0 \ar[r] & A_{p^{m-1}} \ar@{=}[d] \ar[r] & A_{p^m} \ar@{^(->}[d]  \ar[r]^{p^{m-1} \cdot } & A_p \ar@{^(->}[d] \ar[r] & 0 \\
0 \ar[r] & A_{p^{m-1}} \ar[r] & A \ar[r]^{p^{m-1} \cdot } & A \ar[r] & 0 
}
\]
(It is the limitation of assumption (iv) that forces the induction to terminate at $n$). Again the snake lemma applied to 
\eqref{eq:almostSel}, more precisely a diagram chase, yields exactness of 
\[
\rH^1_{\Sel^Q}(k,A_{p^{m-1}})  \to \rH^1_{\Sel^Q}(k,A_{p^m}) \to \rH^1_{\Sel^Q}(k,A_{p}) 
\]
so that with Step 1 we deduce the theorem by induction on $m$.
\end{pro}

\begin{rmk}
It seems difficult to set up an inductive argument to prove $\Sha^1(k,A_{p^m}) = 0$ for all $m \leq n$ directly. 
\end{rmk}

We are now ready to give a proof of Theorem D from the introduction.

\begin{thm} \label{thm:proofofthmD}
Let $A/k$ be an abelian variety over an algebraic number field, and let $p$ be  a prime number. 
Then 
\[
\Sha^1(k,A_{p^n}) = 0
\]
for all $n\geq 0$, if we assume that 
\begin{enumerate}
\item[(i)] $\rH^1(k(A_p)/k,A_p) = 0$, with the splitting field $k(A_p)$ of the $p$-torsion $A_p$ of $A$, and 
\item[(ii)] the $\Gal_k$-modules $A_p$ and $\End(A_p)$ have no common irreducible subquotient.
\end{enumerate}
\end{thm}
\begin{pro}
We are going to show that an auxiliary finite set of finite primes $Q$ exists (depending on $n$) such that conditions (ii)--(iv) of 
Proposition~\ref{prop:structureofgenselmer} hold. Then 
\[
\Sha^1(k,A_{p^n}) \subseteq \rH^1_{\Sel^Q}(k,A_{p^n}) = 0
\]
by Proposition~\ref{prop:structureofgenselmer} proves the theorem.

\smallskip

First we prove that the Selmer splitting field $L$ of $A/k$ and $k(A_{p^n})$ are linearly disjoint over $k(A_p)$. Indeed, let $K = L \cap k(A_{p^n})$ be their intersection and set $\ov{M} = \Gal(K/k(A_p))$ for the abelian Galois group over $k(A_p)$. Then, since $K/k$ is  Galois, the projection
\[
\Gal(L/k(A_p)) \surj \ov{M}
\]
is a surjection of $\Gal(k(A_p)/k)$-modules. It follows from Lemma~\ref{lem:selmerslpittingfield} that $\ov{M}$ has a composition series as $\Gal(k(A_p)/k)$-module consisting of irreducible subquotients of $A_p$.
On the other hand, the group $\Gal(k(A_{p^n})/k(A_p))$ is a subgroup of 
\[
N = \ker\big(\GL(A_{p^n}) \to \GL(A_p)\big).
\]
The group $N$ is solvable with abelian subquotients 
\[
N_m = \ker\big(\GL(A_{p^m}) \to \GL(A_{p^{m-1}})\big)
\]
that are canonically $\Gal(k(A_p)/k)$-modules and isomorphic to the adjoint representation of $\Gal(k(A_p)/k)$ on $\End(A_p)$. 
Since, by assumption, $A_p$ and $\End(A_p)$ share no irreducible subquotient, we deduce that $\ov{M} = 0$ and $K = k(A_p)$, which means that $L$ and $k(A_{p^n})$ are linearly disjoint over $k(A_p)$.

\smallskip

The Chebotarev density theorem enables us to choose a finite set $Q$ of finite places $v \nmid p$ in the locus of good reduction of $A/k$, so that the Frobenius elements $\Frob_v$ for $v \in Q$ satisfy
\begin{enumerate}
\item[(a)] the image of $\Frob_v$ in $\Gal(k(A_{p^n})/k)$ is trivial,
\item[(b)] the  images of $\Frob_v$ for $v \in Q$ generate $\Gal(L/k(A_p))$.
\end{enumerate}
The linear disjointness of $L$ and $k(A_{p^n})$ over $k(A_p)$ implies that (a) and (b) do not contradict each other.
This shows  (ii)--(iv) of Proposition~\ref{prop:structureofgenselmer} for the set $Q$ and concludes the proof.
\end{pro}

\begin{rmk}
Let $A/k$ be an abelian variety over a number field $k$, and let $n \in \bN$. In \cite{dz:start} Dvornicich and Zannier started the investigation, actually for general commutative algebraic groups, of whether an element $a \in A(k)$ that is locally divisible by $n$ in $A(k_v)$ for (almost) all places $v$ of $k$ must necessarily be globally divisible by $n$, i.e, $a \in n A(k)$. The vanishing criterion of Theorem~\ref{thm:proofofthmD} above together with the exact sequence \eqref{eq:locdivpointssequence} gives a criterion for the equality
\begin{equation} \label{eq:locdivpoints}
\{a \in A(k) ; \ a \in p^n \cdot A(k_v) \text{ for all $v$ }\} = p^n \cdot A(k)
\end{equation}
to hold for a prime number $p$ and all $n \geq 1$. Below in Section \S\ref{sec:caseofellipticcurves} we will verify this criterion in various cases for elliptic curves and thus reprove results obtained on cases where \eqref{eq:locdivpoints} holds, namely \cite{paladinoranieriviada} Corollary 2.
\end{rmk}

%-----------------------------------------------------------------------------------------------------------------
\subsection{The criterion in the case of elliptic curves}

Based on the group theory of $\GL_2$ in Section~\S\ref{sec:group} we can show the following criterion for a $p$-primary 
local--to--global principle for divisibility in the special case of elliptic curves. 

\begin{thm} \label{thm:localglobaldivforellipticcurves}
Let $E/k$ be an elliptic curve and let $p$ be a prime number.
We assume that the following holds for $G =  \Gal(k(E_p)/k) \subseteq \GL(E_p)$.
\begin{enumerate}
\item[(i)] The group $G$ is not contained in a subgroup of $\GL(E_p)$ that is isomorphic to the symmetric group $S_3$, in particular $p >2$; and 
 \item[(ii)] If $E_p$ is a reducible $G$-module, then its semisimplification 
 $E_p^{\rm ss} $ is not of the form
 \begin{enumerate}
 \item $\one \oplus \epsilon_p$ with the mod $p$ cyclotomic character $\epsilon_p$, or 
 \item $\chi \oplus \chi^{2}$  for a character $\chi: \Gal_k \to \bF_p^\ast$ such that $\chi^3 = \epsilon_p$.
 \end{enumerate}
 \end{enumerate}
Then the elements of $\rH^1_{\divisor}(k,E)$ and in particular of $\Sha(E/k)$ are $p$-divisible in $\rH^1(k,E)$.
\end{thm}

\begin{pro}
Since the determinant of the mod $p$ representation $\ov{\rho}_{E} : \Gal_k \to \GL(E_p)$ is the cyclotomic character $\epsilon_p$, conditions (i) and (ii) guarantee by 
Theorem~\ref{thm:groupcrit}  that 
\begin{itemize}
\item $E_p$ and $\End(E_p)$ have no common irreducible factor as $\Gal_k$-modules, 
\item and $\rH^1(G,E_p) =0$. 
\end{itemize}
We deduce by Theorem~\ref{thm:proofofthmD} that $\Sha^1(k,E_{p^n}) = 0$ for all $n \geq 0$, and thus 
\[
\rH^1_{\divisor}(k,E)_{p^\infty} \subseteq p^n \cdot \rH^1(k,E)
\]
for all $n \geq 0$ by Proposition~\ref{prop:prdiv}. This concludes the proof.
\end{pro}

%-----------------------------------------------------------------------------------------------------------------
\subsection{Divisibility by almost all primes --- review of Bashmakov's results}   \label{sec:almostall}
Bashmakov gives a partial answer to Cassels' question, at least when $p$ is large. In \cite{bashmakov:survey} \S5, the treatment is said to be restricted to CM-abelian varieties, but is only carried out for elliptic curves.

\begin{prop}[Bashmakov \cite{bashmakov:survey} Proposition 22] \label{prop:bashmakovalmostall}
Let $k$ be an algebraic number field and let $E/k$ be an elliptic curve. Then for almost all $p$ we have 
$\Sha(E/k)_{p^\infty} \subseteq \divisor(\rH^1(k,E))$. In particular, for such $p$, every class in $\Sha(E/k)$ is $p$-divisible in $\rH^1(k,E)$.
\end{prop}
Bashmakov's main tool in the proof of Proposition~\ref{prop:bashmakovalmostall} is provided by a result of Serre \cite{serre:pointsrationnels} that the image of $\Gal_k$ acting on the full Tate module $\rT E = \prod_\ell \rT_\ell E$ contains an open subgroup of the diagonal torus. For an abelian variety instead of an elliptic curve but restricted to the $\ell$-component this is due to Bogomolov. The bound for $p$, such that Proposition~\ref{prop:bashmakovalmostall} holds, provided by Bashmakov's method depends on the index of this open subgroup in the diagonal torus.
In the meantime, Serre has improved his result so that this approach now shows the following.

\begin{thm} \label{thm:abelianalmostall}
Let $k$ be an algebraic number field and let $A/k$ be an abelian variety. Then for almost all $p$ we have $ \rH^1_{\divisor}(k,A)_{p^\infty} = \divisor(\rH^1(k,A))_{p^\infty}$. In particular, for such $p$, every class in $\Sha(A/k)$ is $p$-divisible in $\rH^1(k,A)$.
\end{thm}

\begin{pro} 
Let $\rho_{A^t,p} : \Gal_k \to \GL(\rT_p(A^t))$ be the Galois representation on the $p$-adic Tate module of the dual abelian variety $A^t$. Then Serre, \cite{serre:lettertoribet} \S2, has shown that there is an integer $N \geq 1$ independent of $p$ such that $\rho_{A^t,p}(\Gal_k)$ contains the diagonal $(\bZ_p^\ast)^N$.  Thus, for $p > N+1$, the image $G=G_{1,p}(A^t)$ of the mod $p$ representation $\ov{\rho}_{A^t,p} : \Gal_k \to \GL(A^t_p)$ meets the center of $\GL(A^t_p)$ non-trivially.
By Lemma~\ref{lem:center} we find that for $p > N+1$ the modules $A_p^t$ and $\End(A_p^t)$ have no common irreducible factor and also $\rH^1(G,A^t_p) = 0$. We conclude by Theorem~\ref{thm:proofofthmD} that
\[
\Sha^1(k,A^t_{p^n}) = 0 
\]
for all $n \geq 0$. By Proposition~\ref{prop:prdiv} this completes the proof.
\end{pro}

\begin{rmk}
It is tempting to ask, whether the bound for $p$ in 
Theorem~\ref{thm:abelianalmostall} depends only on $k$ and maybe $\dim(A)$ but not on the particular abelian variety $A/k$. Theorem~\ref{thm:numberfieldEll} (2) provides such a uniform bound on $p$ in terms of only the degree of $k/\bQ$ in the case of elliptic curves.
\end{rmk}

%%%%%%%%%%%%%%%%%%%%%%%%%%%%%%%%%%%%%%%%%%
\section{The case of elliptic curves} \label{sec:caseofellipticcurves}

%----------------------------------------------------------------------------------------------------------------
\subsection{Applications  to elliptic curves over number fields}

We first analyze the conditions asked by  Theorem~\ref{thm:localglobaldivforellipticcurves} in the case of elliptic curves over  arbitrary number fields.

\begin{cor} \label{cor:ell}
Let $k$ be a number field. Then for 'most'  elliptic curves $E/k$, in particular for infinitely many of them, the group  $\rH^1_{\divisor}(k,E)$ and in particular $\Sha(E/k)$ is $p$-divisible in $\rH^1(k,E)$ for all odd primes $p$. 
\end{cor}

\begin{pro}
Jones \cite{Jo10} (for $k=\bQ$) and Zywina \cite{Zy10} consider the set of elliptic curves 
\[
Y^2=X^3+aX+b
\]
such that $a$ and $b$ are integers of $k$ and $h(a,b)< x$ (where $h$ denotes a height on such pairs). They show (\cite{Zy10} Theorem 1.3 and Theorem 1.6) that the ratio of the cardinality of the subset of elliptic curves such that 
\begin{equation} \label{eq:bigG1E}
\frac{|\SL_2(\bF_p)|}{|\Gal(k(E_p)/k) \cap \SL_2(\bF_p)|} \leq 2
\end{equation}
for every prime $p$ by the total number of curves  in the box $h(a,b)< x$ approaches $1$ as $x$ goes to infinity.  When \eqref{eq:bigG1E} holds and $p > 2$, then $\Gal(k(E_p)/k)$ can neither be contained in an $S_3$ nor in a Borel subgroup of $\GL(E_p)$.
Hence in this sense  for 'most'  elliptic curves $E/k$ the conditions of Theorem~\ref{thm:localglobaldivforellipticcurves} hold for every prime $p>2$.
\end{pro}

\smallskip

Recall from Section~\S\ref{sec:piofd} that $\pi(d)$ is the maximal prime number such that an elliptic curve $E$ over a number field $k/\bQ$ of degree $d$ admits a non-trivial $k$-rational $p$-torsion point.

\begin{thm} \label{thm:numberfieldEll}
Let $E/k$ be an elliptic curve defined over an algebraic number field $k$ of degree $d$ over $\bQ$. Then the following holds.
\begin{enumerate}
\item The group
$\rH^1_{\divisor}(k,E)$ and therefore $\Sha(E/k)$ is $p$-divisible in $\rH^1(k,E)$ for a prime number $p \geq 5$ under the following conditions:
\begin{enumerate}
\item[(i)] The extension $k(\zeta_p)/k$ has degree $\geq 3$, and 
\item[(ii)] no elliptic curve $E'$ which is $k$-isogenous to $E$ has a $k$-rational $p$-torsion point, in particular if $p> \pi(d)$, and
\item[(iii)] $p> (Nv + \sqrt{Nv})^2$ where $Nv = |\bF_v|$ is the size of the residue class field $\bF_v$ for a place $v \nmid 3p$ of $k$. 
\end{enumerate}
\item If $p > \max\{(2^d  + 2^{d/2})^2, \pi(d)\}$,
then $\rH^1_{\divisor}(k,E)$ and therefore $\Sha(E/k)$ is $p$-divisible in $\rH^1(k,E)$.
\end{enumerate}
\end{thm}
\begin{pro}
(1) Since the determinant of the mod $p$ representation $\ov{\rho}_{E} : \Gal_k \to \GL(E_p)$ is the cyclotomic character $\epsilon_p$,  property (i) shows that $\det(\Gal(k(E_p)/k))$ has order $\geq 3$. This implies $\Gal(k(E_p)/k)$  is not contained in a subgroup of $\GL(E_p)$ that is isomorphic to the symmetric group $S_3$ (see  Corollary~\ref{cor:detS3}). Hence, we are only concerned by the second condition of 
Theorem~\ref{thm:localglobaldivforellipticcurves}. Since (ii) implies that $E_p^{\rm ss} \not \cong \one \oplus \epsilon_p$, it remains to verify that $E_p^{\rm ss} $  is not of the form  $\chi \oplus \chi^{2}$  for some character $\chi: \Gal_k \to \bF_p^\ast$ such that $\chi^3 = \epsilon_p$.

\smallskip

We argue by contradiction. Note that  the ramification at $v \nmid 3p$ of a character $\chi$ that solves $\chi^3 = \epsilon_p$ is at most tame and of degree $e_v \mid 3$. Let $k'/k$ be an extension of degree $e_v$ depending on $v$ with a place $w \mid v$ of $k'$ such that $k'_w/k_v$ is totally tamely ramified of degree $e_v$. In particular, the size $Nw$ of the residue field $\bF_w$ at $w$ agrees with $Nv$.
It follows then from Abhyankar's Lemma and $E_p^{\rm ss} = \chi \oplus \chi^2$ that the inertia group $\rI_w \subset \Gal_{k'}$ acts unipotently on $E_p = \rT_p(E)/p\rT_p(E)$. Since a priori the action of inertia on $\rT_p(E)$ is quasi-unipotent, Lemma~\ref{lem:unipotentmodp} below applies, and $\rI_w$ must even act unipotently. The criterion of semistable reduction  \cite{sga7-1} IX \S3 Proposition 3.5 then shows that $E/k$ has semistable reduction at $w$ after scalar extension to $k'$. 

If $E/k$ has multiplicative reduction in $w$, then there is an unramified quadratic character $\delta : \Gal_{k'_w} \to \{\pm 1\}$ such that as a $\Gal_{k'_w}$-representation
\[
E_p^{\rm ss}|_{\Gal_{k'_w}} \cong \delta \oplus \delta \epsilon_p \cong  \chi \oplus \chi^2.
\]
If $\delta = \chi$ and $\delta \epsilon_p = \chi^2$, then $\delta = \chi = \epsilon_p$ and thus $\epsilon_p(\Frob_w) = \pm 1 \in \bF_p^\ast$, which means
\[
p \mid Nw \pm 1 = Nv \pm 1.
\]
If on the other hand $\delta = \chi^2$ and $\delta \epsilon_p = \chi$, then $\delta = \chi^2 = \epsilon_p^2$ and thus $\epsilon_p(\Frob_w^2) = \pm 1 \in \bF_p^\ast$, which means
\[
p \mid (Nw)^2 \pm 1 = (Nv)^2 \pm 1.
\]

Finally, if $E/k$ has good reduction in $w$, then $E_p$ is an unramified $\Gal_{k'_w}$-module and 
\[
a = a_w = 1 + Nw - |\ov{E}(\bF_w)| \equiv \tr(\Frob_w| E_p^{\rm ss}) \mod p
\]
Since $E_p^{\rm ss}|_{\Gal_{k'_w}}  \cong  \chi \oplus \chi^2$ where $\chi^3 = \epsilon_p$ (as in Mazur's argument presented in \cite{serre:pointsrationnels}) we deduce that 
\[
p \mid \big(\chi(\Frob_w)^2 + \chi(\Frob_w)\big)^3 - a^3 \equiv Nv + (Nv)^2 + 3Nv \cdot a - a^3 \mod p
\]
with $a \in \bZ$ and $|a| \leq 2\sqrt{Nv}$ due to the Hasse--Weil bound. An analysis of the function $f(x) = 3Nx-x^3$ on the interval $-2\sqrt{N} \leq x \leq 2 \sqrt{N}$ shows that 
\[
|3Nv \cdot a  - a^3| \leq 2 \sqrt{(Nv)^3}.
\]
Consequently, we find
\[
0 < (Nv -\sqrt{Nv})^2 \leq Nv + (Nv)^2 + 3Nv \cdot a - a^3 \leq (Nv + \sqrt{Nv})^2
\]
which leads in the case of good reduction at $w$  to 
\[
p \leq  (Nv + \sqrt{Nv})^2.
\]
Consequently, if $p > (Nv + \sqrt{Nv})^2$ then $E_p^{\rm ss} \not \cong \chi \oplus \chi^ 2$ as claimed.

\smallskip

(2) Observe that if $p > \max\{(2^d  + 2^{d/2})^2, \pi(d)\}$, then $p$ must be odd and hence condition (iii) holds with respect to any place $v \mid 2$.  
Moreover, since $p > \pi(d)$, condition (ii) holds by definition of $\pi(d)$. 
Finally, for (i) we use the crude estimate
\[
[k(\zeta_p):k]  \geq \bruch{[\bQ(\zeta_p):\bQ]}{[k:\bQ]} = \bruch{p-1}{d} > \bruch{(2^d  + 2^{d/2})^2 -1}{d} 
> \bruch{4^d}{d} \geq \bruch{1 + d \cdot 3}{d} > 3.
\]
Therefore, (2) holds as a special case of (1).
\end{pro}

\begin{lem} \label{lem:unipotentmodp}
Let $p$ be an odd prime number. Let $\rho : G \to \GL_2(\bZ_p)$ be a quasi-unipotent continuous representation of a pro-finite group $G$ such that the mod $p$ representation $\ov{\rho} : G \to \GL_2(\bF_p)$ is unipotent. If $\rho$ is not unipotent, then $p =3$ and $\rho(G) \cong \bZ/3\bZ$.
\end{lem}
\begin{pro}
Since we assume that $G$ acts quasi-unipotently, there is an open normal subgroup $G^0 \subset G$ such that $\rho|_{G^0}$ is unipotent. If $G^0$ acts non-trivially, then it fixes a unique $\bZ_p$-line, and $\rho(G)$ lies in a Borel subgroup, namely the normalizer of the stabilizer of the line.  Due to $\rho$ being quasi-unipotent and $\ov{\rho}$ being unipotent, the corresponding diagonal characters map to a finite group of  $1$-units in $\bZ_p^\ast$ and therefore are trivial. Hence the representation $\rho$ is unipotent. 

If on the other hand $\rho(G^0) = 1$, then $\rho(G)$ is a torsion subgroup of $\GL_2(\bZ_p)$ and therefore isomorphic to its image in $\GL_2(\bF_p)$ (see Lemma 9 in \cite{soule:arith}). Since we may assume that $\rho$ is non-trivial we conclude that $\rho(G) \cong \bZ/p\bZ$, and that we have a non-trivial $M \in \GL_2(\bQ_p)$ with $M^p = 1$. The characteristic polynomial of $M$ is quadratic over $\bQ_p$ and its roots are non-trivial $p$th roots of unity. Hence the $p$th cyclotomic extension $\bQ_p(\zeta_p)$ is at most quadratic over $\bQ_p$, whence $p\leq 3$.
\end{pro}

\begin{cor} \label{cor:serretrickEll}
Let $E/k$ be an elliptic curve defined over an algebraic number field $k$, and let $p$ be a prime number. Then $\rH^1_{\divisor}(k,E)$ is $p$-divisible in $\rH^1(k,E)$ under the following conditions:
\begin{enumerate}
\item The extension $k(\zeta_p)/k$ has degree $\geq 3$, and
\item one of the following holds:
\begin{enumerate}
\item $E$ has good reduction above a place $v$ of $k$ with norm $Nv = 3$ and $p > 11$.
\item no elliptic curve $E'$ which is $k$-isogenous to $E$ has a $k$-rational $p$-torsion point and $3 \mid p-1$ but the degree of $\bQ(\zeta_p) \cap k$ over $\bQ$ is prime to $3$. 
\item $3 \mid p-1$, the prime $p$ is unramified in $k/\bQ$, and $p > \pi(d)$ where $d$ is the degree of $k/\bQ$. 
\end{enumerate}
\end{enumerate}
\end{cor}

\begin{rmk}
The conditions (b) and (c) in Corollary \ref{cor:serretrickEll} both already imply condition (1).
\end{rmk}

\begin{pro}[Proof of Corollary \ref{cor:serretrickEll}]
As in the proof of Theorem~\ref{thm:numberfieldEll} we deduce from (1) that the first condition of 
Theorem~\ref{thm:localglobaldivforellipticcurves} holds. Hence, we are only concerned by the second condition of 
Theorem~\ref{thm:localglobaldivforellipticcurves}. 
We therefore assume that $E_p$ is a reducible $\Gal_k$-representation with semi-simplification $E_p^{\rm ss} = \chi_1 \oplus \chi_2$. 
In order to satisfy the second condition of 
Theorem~\ref{thm:localglobaldivforellipticcurves} we have to exclude the following two cases:
\begin{enumerate}
\item[(A)]  $E_p^{\rm ss} = \one \oplus \epsilon_p$,
\item[(B)] $E_p^{\rm ss} = \chi \oplus \chi^2$ with  $\chi^3 = \epsilon_p$.
\end{enumerate}
Under the assumption (c) we have $\bQ(\zeta_p) \cap k = \bQ$ and the assumption in (b) on $k$-rational torsion points holds again by definition of $\pi(d)$.
So (c) in fact implies (b). Under the assumption (b) the option (A) is excluded and the equation $\chi^3 =\epsilon_p$ has no solution, which excludes option (B).

It remains to discuss assumptions (a). Let $v$ be a place of $k$ with residue field $\bF_v$ of cardinality $Nv$ 
and where $E$ has good reduction $\ov{E}/\bF_v$, then as in the proof of Theorem~\ref{thm:numberfieldEll} 
we have $a \in \bZ$ with $|a| \leq 2\sqrt{Nv}$ and 
\[
a = \tr(\Frob_v| \rT_p(E)) \equiv \chi_1(\Frob_v) + \chi_2(\Frob_v) \mod p.
\]
We distinguish the two cases and argue similarly to the proof of Theorem~\ref{thm:numberfieldEll}.
\begin{enumerate}
\item[(A)]  If  $E_p^{\rm ss} = \one \oplus \epsilon_p$, then 
\[
1+ Nv - a \equiv 0 \mod p.
\]
\item[(B)] If $E_p^{\rm ss} = \chi \oplus \chi^2$ with  $\chi^3 = \epsilon_p$, then 
\[
 (Nv +(Nv)^2) + 3 a \cdot Nv - a^3 \equiv  \big(\chi(\Frob_v)^2 + \chi(\Frob_v)\big)^3 - a^3 \equiv   0 \mod p.
\]
\end{enumerate}
We finally exploit assumption (a), so that $Nv = 3$ and $|a| \leq 3$. The constraints in the two cases now are
\begin{enumerate}
\item[(A)]  $p$ divides one of $4-a \in\{1,2,3,4,5,6,7\}$,
\item[(B)]  $p$ divides one of $12 + 9a - a^3 \in \{12, 2, 4,12,20,22,12\}$.
\end{enumerate}
Hence, if $p>11$ none of the cases (A) or (B) can occur, and thus $E/k$ satisfies the assumptions of 
Theorem~\ref{thm:localglobaldivforellipticcurves} for $p$.
\end{pro}

\begin{cor} \label{cor:irreduciblecase}
Let $E/k$ be an elliptic curve defined over an algebraic number field $k$, and let $p$ be an odd prime number. Then  
$\rH^1_{\divisor}(k,E)$ is $p$-divisible in $\rH^1(k,E)$ if $E_p$ is an irreducible $\Gal_k$-representation and $[k(\zeta_p):k] \not=2$.
\end{cor}
\begin{pro}
We have to verify the conditions of Theorem~\ref{thm:localglobaldivforellipticcurves}. Since $E_p$ is assumed irreducible, condition (ii) is automatic. In order to satisfy also condition (i) of Theorem~\ref{thm:localglobaldivforellipticcurves} we argue by contradiction. If $\Gal(k(E_p)/k) \subseteq \GL(E_p)$ lies in a copy of $S_3$, then for $E_p$ to be irreducible we must have $p \geq 5$ and $\Gal(k(E_p)/k) \cong  S_3  \subseteq \GL_2(\bF_p)$. But then 
\[
[k(\zeta_p):k] = |\det(\Gal(k(E_p)/k)| = |\det(S_3)| = 2
\]
by Corollary~\ref{cor:detS3} and this was excluded. 
\end{pro}

%----------------------------------------------------------------------------------------------------------------
\subsection{Quadratic twists}

For an elliptic curve $E/k$  and a character 
\[
\tau : \Gal_k \to \{\pm 1\} \subseteq \Aut(E)
\]
we denote by $E^\tau$ the quadratic twist of $E$ by $\tau$. The  $p$-torsion of the twist  $E^\tau$ is given by $(E^\tau)_p = E_p \otimes \tau$ as a $\Gal_k$-module.

\begin{cor} \label{cor:quadratictwist}
Let $k$ be a number field and let $E/k$ be an elliptic curve. 
\begin{enumerate}
\item Let $p \geq 3$ be a fixed prime number. Then among the quadratic twists of $E$ there are at most $3$ twists $E^\tau$ such that $\rH^1_{\divisor}(k,E^\tau)$ is not $p$-divisible in $\rH^1(k,E^\tau)$.
\item For all but a finite number of quadratic twists $E^\tau$ of $E$ the group $\rH^1_{\divisor}(k,E^\tau)$ is $p$-divisible in $\rH^1(k,E^\tau)$ for all $p \geq 3$.
\end{enumerate}
\end{cor}
\begin{pro}
(1) We verify the conditions of Theorem~\ref{thm:localglobaldivforellipticcurves} for almost all quadratic twists.  Since twisting is transitive, we may first assume that $G = \Gal(k(E_p)/k) \subseteq \GL(E_p)$ is isomorphic to one of the groups
\begin{enumerate}
\item[(a)] the trivial group,
\item[(b)] $\bZ/2\bZ$,  non-centrally in $\GL(E_p)$ (see Corollary \ref{cor:detS3}), 
\item[(c)] $\bZ/3\bZ$,
\item[(d)] $S_3$. 
\end{enumerate}
For a non-trivial $\tau : \Gal_k \to \{\pm 1\}$ we find $G^\tau = \Gal(k(E^\tau_p)/k) \subseteq \GL(E^\tau_p)$ equal to respectively
\begin{enumerate}
\item[(a)] $\bZ/2\bZ$, centrally in $\GL(E_p^\tau)$,
\item[(b)] $\bZ/2\bZ \times \bZ/2\bZ$,
\item[(c)] $\bZ/3\bZ \times \bZ/2\bZ$,
\item[(d)] $S_3 \times \bZ/2\bZ$,  since a transposition in $S_3$ cannot be central in $\GL(E^\tau_p)$.
\end{enumerate}
Consequently, every non-trivial twist $E^\tau$ satisfies condition (1) of Theorem~\ref{thm:localglobaldivforellipticcurves}, which means that in general, at most one of the twists of $E/k$ can fail condition (1) of Theorem~\ref{thm:localglobaldivforellipticcurves}.

Secondly, we assume that $E_p$ is reducible and that condition (2) of Theorem~\ref{thm:localglobaldivforellipticcurves} fails. If a non-trivial twist $E^\tau$ also fails condition (2), then twisting by $\tau$ either preserves $\one \oplus \epsilon_p$ or $\chi \oplus \chi^2$ with $\chi^3 = \epsilon_p$, or it transforms one to the other. In any case, the character $\tau \not= \one $ is one of the following list 
\begin{enumerate}
\item[(a)]  $(\one \oplus \epsilon_p) \otimes \tau \cong \one \oplus \epsilon_p$: then $\tau = \epsilon_p$,
\item[(b)]  $(\chi \oplus \chi^2) \otimes \tau \cong  \chi \oplus \chi^2$: then $\tau = \chi= \epsilon_p$,
\item[(c)] $(\one \oplus \epsilon_p) \otimes \tau \cong \chi \oplus \chi^2$: then $\tau = \chi = \epsilon_p$
or $\tau = \chi^2 = \epsilon_p^2$,
\item[(d)] $(\chi \oplus \chi^2) \otimes \tau \cong \one \oplus \epsilon_p$: then $\tau = \chi = \epsilon_p$ or $\tau = \chi^2 = \epsilon_p^2$.
\end{enumerate}
In any case $\tau = \epsilon_p$ or $\epsilon_p^2$ which cannot  both be non-trivial quadratic characters. Hence, condition (2) of Theorem~\ref{thm:localglobaldivforellipticcurves} can in general only fail for at most $2$ twists of $E/k$.

\smallskip

(2) Due to (1) it suffices to show that for $p \gg 3$ and all quadratic twists $E^\tau/k$ the group $\rH^1_{\divisor}(k,E^\tau)$ is $p$-divisible in $\rH^1(k,E^\tau)$.  If $E/k$ has no CM, then by \cite{Serre}  \S4.4 we have 
\[
\Gal(k(E_p)/k) = \GL(E_p)
\]
for $p \gg 3$. If $E/k$ has CM, then by  \cite{Serre} \S4.5 the image $\Gal(k(E_p)/k) \subseteq \GL(E_p)$ contains a full (split or non-split) torus $C \subset \GL(E_p)$ for $p \gg 3$. In both cases $\Gal(k(E_p)/k)$ contains the group of diagonal matrices $\cong \bF_p^\ast$ and continues to contain at least the squares $\cong (\bF_p^\ast)^2 \not=1$ after twisting. Lemma~\ref{lem:center} then provides that the conditions of 
Theorem~\ref{thm:localglobaldivforellipticcurves} are satisfied for all quadratic twists of $E$.
\end{pro}

\begin{cor} \label{cor:allj}
Let $k$ be a number field. For every $j \in k$ there is an elliptic curve $E/k$ with $j$-invariant $j$ such that $\rH^1_{\divisor}(k,E)$ and therefore $\Sha(E/k)$ is $p$-divisible in $\rH^1(k,E)$ for all odd primes $p$.
\end{cor}
\begin{pro}
Quadratic twists share the same $j$-invariant, and Corollary~\ref{cor:quadratictwist}.
\end{pro}

\begin{rmk}
Corollary~\ref{cor:allj} would be automatic if among the quadratic twists of a given $E/k$ we could find a curve $E^\tau$ with $\Sha(E^\tau/k) = 0$ or at least of order a power of  $2$.  
\end{rmk}

%----------------------------------------------------------------------------------------------------------------
\subsection{Applications  to elliptic curves over \texorpdfstring{$\bQ$}{the rationals}} \label{sec:elloverQ}

Now we attempt to find optimal results in the case of $k=\bQ$.

\begin{thm} \label{thm:rationalEll}
Let $E/\bQ$ be an elliptic curve defined over the rationals. Then the following holds.
\begin{enumerate}
\item  $\rH^1_{\divisor}(\bQ,E)$ and therefore $\Sha(E/\bQ)$ is $p$-divisible in $\rH^1(\bQ,E)$ for all primes $p>7$.
\item \label{item:badpairsEp}
Let $p$ be an odd prime number such that $\rH^1_{\divisor}(\bQ,E)$ is not $p$-divisible in $\rH^1(\bQ,E)$. Then we have one of the following cases.
\begin{enumerate}
\item[(a)] $p=3$ and $E_3^{\rm ss} = \one \oplus \epsilon_3$, 
\item[(b)] $p=5$ and $E_5^{\rm ss} = \one \oplus \epsilon_5$, or $E_5^{\rm ss} = \epsilon_5^3 \oplus \epsilon_5^2$,  
\item[(c)] $p=7$ and $E_7^{\rm ss} = \one \oplus \epsilon_7$.
\end{enumerate}
In particular, in all the above cases the inertia group $\rI_v \subset \Gal_\bQ$ for a place $v \nmid p$ acts through a $p$-group on $\rT_p(E)$. Moreover, in case (b) and (c) the curve $E$ has semistable reduction outside of $p$.
\item  $\rH^1_{\divisor}(\bQ,E)$ and therefore $\Sha(E/\bQ)$ is $p$-divisible in $\rH^1(\bQ,E)$ for all  odd primes $p$ where $E$ has supersingular or non-split multiplicative reduction at $p$.  
\end{enumerate}
\end{thm}
\begin{pro}
(1) Observe that by Theorem \ref{thm:numberfieldEll} (2) we know that $\rH^1_{\divisor}(\bQ,E)$ and hence $\Sha(E/\bQ)$ are $p$-divisible in $\rH^1(\bQ,E)$ for all primes $p > \max\{(2+\sqrt 2)^2,\pi(1)\}$. Since $11<(2+\sqrt 2)^2<12$ and $\pi(1) = 7$ by Mazur \cite{mazur:ratisog} Theorem~2, it follows that we can now restrict our attention to the case $p=11$.

Again by work of Mazur (see Theorem~2 in \cite{mazur:ratisog}) we know that $E(\bQ)_{11}$ is trivial and hence all we need is to identify (and deal differently with) all $E/\bQ$ such that  
$E_{11}^{\rm ss}$ is of the form $\chi \oplus \chi^{2}$  for a character $\chi: \Gal_k \to \bF_{11}^\ast$ such that $\chi^3 = \epsilon_{11}$. Observe that  $\chi^3 = \epsilon_{11}$ implies that $\chi= \epsilon_{11}^{\otimes 7}$. 
We know that elliptic curves such that $E_{11}$ is reducible correspond up to quadratic twist to non-cuspidal rational points\footnote{We thank Brian Conrad for suggesting this way of dealing with the prime $11$.} of $X_0(11)$ of which there are three (see \cite{modularfunctionsIV} page 79). These three rational points correspond up to quadratic twist to the following three elliptic curves $E/\bQ$ of conductor $121$ (the code is as in Cremona's list \cite{cremona:data}): 

\begin{center}
\begin{tabular}[t]{|c||c|c||c|c|} \hline
&& & &  \\[-2ex]
   elliptic curve & $[a_1,a_2,a_3,a_4,a_6]$ & $j$-invariant  & $\tr(\Frob_2|\rT_{11}(E))$ & $E_{11}^{\rm ss}$ \\[0.5ex] \hline
    121b1 & $[0,-1,1,-7,10]$ & $- 2^{15}$ &  $0$ &  $\epsilon_{11}^{\otimes 3} \oplus \epsilon_{11}^{\otimes 8} $ \\[0.5ex]
    121c1 & $[1,1,0,-2,-7]$ & $- 11^2$ & $1$ &  $\epsilon_{11}^{\otimes 4} \oplus \epsilon_{11}^{\otimes 7} $ \\[0.5ex]
    121c2 & $[1,1,0,-3632,82757]$ & $- 11 \cdot 131^3$ & $1$ & $\epsilon_{11}^{\otimes 4} \oplus \epsilon_{11}^{\otimes 7} $  \\[0.5ex]   \hline
  \end{tabular}
\end{center}

\medskip

\noindent These elliptic curves $E/\bQ$ have good reduction outside $11$.
By reducing the affine equation of the respective elliptic curve modulo $2$ and counting $|\ov E(\bF_2)|$ we compute 
\[
\tr(\Frob_2|\rT_{11}(E)) = 2 + 1 - |\ov E(\bF_2)|
\]
recorded in the table above. On the other hand, since  for these curves $E_{11}$ is reducible and unramified away from $11$ class field theory tells us that $E_{11}^{\rm ss} \cong \epsilon_{11}^a \oplus \epsilon_{11}^b$ for some $a,b \in \bZ/10\bZ$ with $a+b \equiv 1 \mod 10$.
To complete the above table we determine the pair $(a,b)$ by comparing $\tr(\Frob_2|\rT_{11}(E))$ with $2^a + 2^b \mod 11$ as follows.

\begin{center}
\begin{tabular}[t]{|c||c|c|c|c|c|c|} \hline
&& & & & & \\[-2ex]
 $(a,b)$                 & $(0,11)$ & $(1,10)$ & $(2,9)$ & $(3,8)$ & $(4,7)$ & $(5,6)$ \\[0.5ex] \hline
$2^a+2^b \mod 11$ & 3 & 3 & -1 & 0 & 1 & -3 \\[0.5ex]  \hline
  \end{tabular}
\end{center}

\medskip

\noindent 
 If a quadratic twist  $E^\tau/\bQ$ has $E^{\tau, {\rm ss}}_{11} = \epsilon_{11}^{\otimes 7} \oplus \epsilon_{11}^{\otimes 4} \cong E_{11}^{\rm ss} \otimes \tau$, then $\tau$ must be a power of $\epsilon_{11}$, namely $\tau=\one$ or $\tau=\epsilon_{11}^{\otimes 5}$, and so this does not occur for $\tau = \epsilon_{11}^{\otimes 5}$ in view of the above determined structure of  $E_{11}^{\rm ss}$ in the three cases. 
 
Consequently, in the case of $p=11$ we are left with exactly two potential exceptions, the two 
$11$-isogenous non-CM curves labeled ``121c1'' and ``121c2''. For these two curves only, the criterion of 
Theorem~\ref{thm:localglobaldivforellipticcurves} does not apply. 
Let $E$ be one of the elliptic curves labeled ``121c1'' or ``121c2''. Since $E$ has good reduction outside the regular prime $p=11$, and since the Galois action on $E_{11}^{\rm ss}$ factors over $\epsilon_{11}$ we deduce from \cite{ciperianistix:bashmakov}
Proposition~5 with $k=\bQ$ that our potential exceptions obey the theorem as well.

\smallskip

(2) The list follows immediately from Theorem~\ref{thm:localglobaldivforellipticcurves} and the analysis of its conditions in the proof of part (1). Note that $\chi^3 = \epsilon_7$ has no solution. 

We concude that for $v \nmid p$ the inertia group $\rI_v$ acts unipotently on $E_p = \rT_p(E)/p\rT_p(E)$. The action of $\rI_v$ on $\rT_p(E)$ is a priori quasi-unipotent, so that Lemma~\ref{lem:unipotentmodp} applies and the action is in fact unipotent for $p = 5$ or $p=7$. We conclude by the criterion for semistable reduction, see  \cite{sga7-1} IX \S3 Proposition 3.5.

\smallskip

(3) We now consider the case when $E$ has either supersingular or non-split multiplicative reduction at $p$. In this case the second condition of Theorem~\ref{thm:localglobaldivforellipticcurves} holds because 
\begin{enumerate}
\renewcommand{\labelenumi}{(\theenumi)}
\renewcommand{\theenumi}{\roman{enumi}}
\item if $E$ has supersingular reduction at $p$ then $\rho_{E,1}$ is irreducible (see Proposition 2.11 in  \cite{DDR}); and
\item  if $E$ has non-split multiplicative reduction at $p$,  we have that $E_p^{\rm ss}| _{\Gal_{\bQ_p}}= \delta \oplus \delta \epsilon_p$ where $\delta$ is the unique unramified quadratic character of $\Gal_{\bQ_p}$ (see Proposition 2.12 in  \cite{DDR}).
\end{enumerate}
It remains to verify the first condition of  Theorem~\ref{thm:localglobaldivforellipticcurves} for $p=3$. We argue by contradiction and assume that the image of $\Gal_\bQ$ in $\GL(E_3)$ is contained in a subgroup isomorphic to $S_3$. Then the proof of Lemma~\ref{lem:S3} shows that $E_3$ is reducible with semisimplification $E_3^{\rm ss} = \one \oplus \epsilon_3$ and we are back in the case of the second condition of Theorem~\ref{thm:localglobaldivforellipticcurves} that we already dealt with.
\end{pro}

\begin{cor} \label{cor:quadratictwistoverQ}
Among the quadratic twists $E^\tau/\bQ$ of an elliptic curve $E/\bQ$ we find at most  one odd prime number $p=3,5$ or $7$ and at most 
\begin{enumerate}
\item[(a)] $2$ twists for $p=3$,
\item[(b)] $2$ twists for $p=5$,
\item[(c)] $1$ twist for $p=7$,
\end{enumerate}
such that $\rH^1_{\divisor}(\bQ,E^\tau)$ is not $p$-divisible in $\rH^1(\bQ,E^\tau)$. In case (b) and (c) such a twist has semistable reduction at all $\ell \not= p$. In particular, for all but at most $2$ of the quadratic  twists of $E/\bQ$ we have that $\rH^1_{\divisor}(\bQ,E^\tau)$ and therefore $\Sha(E^\tau/\bQ)$ is $p$-divisible in $\rH^1(\bQ,E^\tau)$ for all odd prime numbers $p$.
\end{cor}
\begin{pro}
The proof is straight forward following the proof of part (1) of Corollary~\ref{cor:quadratictwist} together with the knowledge of the precise structure of bad pairs $(E,p)$ for $k=\bQ$ from Theorem ~\ref{thm:rationalEll} \eqref{item:badpairsEp}. It remains to exclude that among the twists of a given $E/\bQ$ more than one of the cases (a)--(c) can occur. Cases (b) and (c) cannot both among the twists for a fixed $E/\bQ$ since no elliptic curve has rational isogeny of degree $35$ by \cite{kenku:rationalisogenies} Theorem 1.

We argue by contradiction. Let us assume that the primes $3$ and $p \in \{5,7\}$ occur as bad primes among the twists of $E/\bQ$. Note that for the prime $p=5$ the two possible bad semisimplifications are twists of each other by $\epsilon_5^2$.  Therefore, after twisting, we may assume that $E_p^{\rm ss} = \one \oplus \epsilon_p$ and that there is a quadratic character $\tau$ with $E_3^{\tau,{\rm ss}} = \one \oplus \epsilon_3$. 

It follows from Lemma \ref{lem:unipotentmodp} 
that an inertia group $\rI_\ell \subset \Gal_{\bQ}$ for $\ell \not=3$ acts on $\rT_3(E^\tau)$
\begin{enumerate}
\item either unipotently  and  $E^\tau$ has semistable reduction at $\ell$,
\item or via a finite $3$-group and $E^\tau$ has potentially good reduction at $\ell$, more precisely with the image of inertia $\Phi_\ell \cong \bZ/3\bZ$.
\end{enumerate}
In the potentially good reduction case and for $\ell \not=p$ we find the same image  $\Phi_\ell$ in $\GL(E^\tau_p)$, namely as the image of inertia in the geometric automorphism group of the special fibre of the potential good reduction of $E^\tau$.  From $E_p^{\tau,{\rm ss}} = \tau \oplus \tau \epsilon_p$ we conclude that the mod $p$ representation $E_p^\tau$ does not allow $I_\ell$ to act via a $3$-group. Thus  $E^\tau/\bQ$ is in fact semistable at $\ell$. Consequently, the Galois module structure of $E_p^{\tau,{\rm ss}}$ now implies that  $\tau$ must be unramified outside $3p$.

If $\tau$ ramifies at $p$, then $E_p^{\tau,{\rm ss}} = \tau \oplus \tau \epsilon_p$ forbids semistable reduction at $p$ for $E^\tau$, whence $E_3^{\tau, {\rm ss}}  = \one \oplus \epsilon_3$ shows that $E^\tau$ has potentially good reduction at $p$. More precisely, good reduction at $p$ occurs after a Galois extension $k/\bQ$ with ramification degree $3$ above $p$ (for example the extension $k = \bQ(E_3^\tau)$ works). But then $E_p^{\tau,{\rm ss}} = \tau \oplus \tau \epsilon_p$  still forbids good ordinary reduction. In case of supersingular reduction, the image of inertia $\rI_{k,p} = \rI_p \cap \Gal_k$ at $p$ (after the extension) is contained in the intersection of a non-split torus in $\GL(E_p)$ with the split torus associated to the decomposition $E_p^{\tau,{\rm ss}} = \tau \oplus \tau \epsilon_p$. This intersection is contained in the central diagonal torus, and therefore $\tau$ and $\tau \epsilon_p$ agree on $\rI_{k,p}$. This leads to $\epsilon_p(\rI_{k,p}) =1$ in contradiction to $\epsilon_p(\rI_p) = \bF_p^\ast$ and $|\rI_p/\rI_{k,p}| = 3$.

It follows that $\tau$ must be unramified outside $3$, and then $\tau = \epsilon_3$ or $\tau = \one$ so that in any case $E_3^{\rm ss} = (\one \oplus \epsilon_3) \otimes \tau = \epsilon_3 \oplus \one$. Hence, there is an elliptic curve $E'/\bQ$ with a $\bQ$-rational point of order $3p$ that is $\bQ$-isogenous to $E$. This  contradicts Mazur's list \cite{mazur:eisenstein} Theorem 8 of possible orders of  rational torsion points.
\end{pro}

\begin{rmk}
(1) It is interesting to note that the difficult odd prime numbers with respect to showing $p$-divisibility for elliptic curves  $E$ over $\bQ$ are exactly the odd \textit{Mazur prime numbers}, i.e., those prime numbers for which $E$ may contain $\bQ$-rational $p$-torsion elements. To some extent this is a consequence of our method, but it is tempting to look for a deeper connection.  

(2) In this respect it is amusing that in order to produce  $p$-torsion in $\Sha(E^\tau/\bQ)$ for quadratic twists of $E/\bQ$, a frequent assumption requires $E$ to have a $\bQ$-rational $p$-torsion point. For example, in \cite{balog-ono} Theorem 2 Balog and Ono show that if $p = 3,5$ or $7$ such that $E$ is good\footnote{Apparently a mild constraint, see \cite{balog-ono}.} with respect to $p$ and has a $\bQ$-rational $p$-torsion point, then with $\tau_d$ ranging over quadratic characters of fundamental discriminant $d$, and $E^d= E^{\tau_d}$, we have an asymptotic lower bound 
\[
\left| \{0 < -d < x \ ; \  E^{d} \text{ has analytic rank } 0 , \text{ and  }  \Sha(E^{d}/\bQ)_p \not= 0 \} \right| \gg \bruch{x^{1/2+1/(2p)}}{\log^2 x}
\]
as $x$ goes to infinity. In particular, an assumption that we like to avoid, namely having rational $p$-torsion, actually helps us to show that our $p$-divisibility result 
Corollary~\ref{cor:quadratictwistoverQ} is non-trivial because the Tate--Shafarevich group in question actually sometimes has non-trivial $p$-torsion.

(3) Matsuno in \cite{matsuno:sha13} Theorem 5.1 shows that $\dim_{\bF_p} \Sha(E/\bQ)_p$ is unbounded for $p=2,3,5,7$ or $13$ among all elliptic curves $E/\bQ$. More precisely, this result also searches for $p$-torsion elements in the Tate--Shafarevich group among quadratic twists of a given elliptic curve with a reducible mod $p$ representation. The list of primes is exactly the list of primes $p$ such that the modular curve $X_0(p)$ is rational.
\end{rmk}

\begin{cor} \label{cor:divisible-allprimes}
Let $E/\bQ$ be an elliptic curve which does not have a rational 4-cyclic subgroup. Then among the quadratic twists of $E$ there are infinitely many $E'/\bQ$ with $\Sha(E'/\bQ)$ divisible in $\rH^1(\bQ,E')$.
\end{cor}

\begin{pro}
Observe that if $E$ has a $4$-cyclic subgroup defined over $\bQ$ then it has a $2$-cyclic subgroup defined over $\bQ$, and thus a non-trivial $\bQ$-rational $2$-torsion point and $\#E_2(\bQ) \geq 2$. Hence we have to consider the following three cases.
\begin{enumerate}
\item$E_2(\bQ)=0$: by Theorem 1.4 of Mazur and Rubin \cite{MR}  we know that $E$ has infinitely many twists $E'$ such that $\rH_\Sel^1(\bQ,E'_2)$ is trivial and consequently so is $\Sha(E'/\bQ)_2$.
\item $E_2(\bQ)\simeq\bZ/2\bZ$ and $E$ has no $4$-cyclic subgroup defined over $\bQ$: 
\begin{enumerate}
\item if $E$ has no $4$-cyclic subgroup defined over $\bQ(E_2)$ then $E$ has infinitely many twists $E'$ such that the $\bZ/2\bZ$-rank of $\rH_\Sel^1(\bQ,E'_2)$ equals $1$ (see Theorem 1.3 of Klagsbrun \cite{Klagsbrun});   
\item if $E$ has a $4$-cyclic subgroup defined over $\bQ(E_2)$ then since elliptic curves over $\bQ$ do not have constant 2-Selmer parity \cite{MR}, by Theorem 1.5 of Klagsbrun \cite{Klagsbrun} we have that $E$ has infinitely many twists $E'$ such that the $\bZ/2\bZ$-rank of $\rH_\Sel^1(\bQ,E'_2)$ equals $1$.  
\end{enumerate} 
Consequently,  $E$ has infinitely many quadratic twists $E'$ such that $\Sha(E'/\bQ)_2$ is trivial.
\item $E_2(\bQ)\simeq(\bZ/2\bZ)^2$ and $E$ has no $4$-cyclic subgroup defined over $\bQ$: by work of Heath-Brown \cite{HB}, Swinnerton-Dyer \cite{SD} and Kane \cite{Kane} (see Theorem 2), we know that $E$ has infinitely many twists $E'$ such that the $\bZ/2\bZ$-rank of $\rH_\Sel^1(\bQ,E'_2)$ equals 2 and hence $\Sha(E'/\bQ)_2$ is trivial.
\end{enumerate}

Then by  Corollary~\ref{cor:quadratictwistoverQ} for an infinite subset of these twists $E'$, in fact all but at most $2$, we have that $\Sha(E'/\bQ)$ is $p$-divisible in $\rH^1(\bQ,E')$ for all primes $p$.
\end{pro}

\begin{cor} \label{cor:legendre}
Let $E/\bQ$ be an elliptic curve with $\bQ$-rational $2$-torsion, i.e., an elliptic curve in twisted  Legendre form $aY^2 = X(X-1)(X-\lambda)$. Then $\rH^1_{\divisor}(\bQ,E)$ and therefore $\Sha(E/\bQ)$ is $p$-divisible in $\rH^1(\bQ,E)$ for all $p\geq 5$.
\end{cor}
\begin{pro}
A quadratic twist of $E/\bQ$ still has $\bQ$-rational $2$-torsion. By Mazur \cite{mazur:eisenstein} Theorem 8, the only $\bQ$-rational torsion of odd order that can occur for $E$ is of order $3$. Hence the cases (b) and (c) of Theorem~\ref{thm:rationalEll} \eqref{item:badpairsEp} cannot occur. Note that upon twisting by the quadratic $\epsilon_5^2$ we interchange the two cases in (b).
\end{pro}

%-----------------------------------------------------------------------------------------------------------------------------
\subsection{An example: the Jacobian of the Selmer curve} \label{sec:selmerexample}
The plane cubic 
\[
S = \{3X^3 + 4Y^3 + 5Z^3 = 0\}
\]
describes Selmer's curve of genus $1$ violating the Hasse principle. Its Jacobian $E=\Pic_S^0$ is an elliptic curve over $\bQ$ of analytic rank $0$ given by the homogeneous equation
\begin{equation} \label{eq:selmercurve}
X^3 + Y^3 + 60Z^3 = 0
\end{equation}
with $[1:-1:0]$ as its origin.
The curve $E$ has  Mordell-Weil group $E(\bQ)=0$, see \cite{cassels:lectures}  \S18 Lemma 2, and $3$-torsion in an exact sequence
\begin{equation} \label{eq:3torsofjacofselmer}
0 \to \mu_3 \to E_3 \to \bZ/3\bZ \to 0,
\end{equation}
which splits over a field $k/\bQ$ if and only if $60$ is a cube in $k$. Moreover, $E$ has complex multiplication by $\fo = \bZ[\zeta_3]$, and the complex multiplication is defined over $\bQ(\zeta_3)$. The curve $S$ and $E$ have good reduction over $U = \Spec(\bZ[1/30])$.

\smallskip

The curve $S$, as a principal homogeneous space under $E$ describes a non-trivial $3$-torsion element $[S] \in \Sha(E/\bQ)$, see \cite{mazur:localtoglobal} I \S4 and \S9. 
Mazur and Rubin determine 
\[
\Sha(E/\bQ) \cong \bZ/3\bZ \times \bZ/3\bZ
\]
unconditionally,  see \cite{mazur:localtoglobal} Theorem~1 and \S9. The Selmer curve $S$ served in \cite{ciperianistix:bashmakov} Theorem~B and \S8 as an example generating a non-trivial intersection
\[
\langle [S] \rangle = \Sha(E/\bQ) \cap \Div(\rH^1(\bQ,E)).
\]
Here we answer the divisibility question of Cassels for $E$ and thus justify \cite{ciperianistix:bashmakov} Remark 25.

\begin{prop} \label{prop:casselsforjacofselmer}
For the Jacobian $E/\bQ$ of the Selmer curve we have
\[
\Sha(E/\bQ) + \Div(\rH^1(\bQ,E))_{3^\infty} = \divisor(\rH^1(\bQ,E))_{3^\infty} = \rH^1_{\divisor}(\bQ,E)_{3^\infty}.
\]
\end{prop}
 
The proof of Proposition~\ref{prop:casselsforjacofselmer} requires some preparation. On $E_{3^\infty}$, the $3$-primary torsion, $\Gal_{\bQ(\zeta_3)}$ acts via $\fo$-linear automorphisms, more precisely $R$-linear automorphisms for 
\[
R = \fo \otimes \bZ_3 = \bZ_3[\zeta_3].
\]
The ring $R$ is a discrete valuation ring, and as  an $R$-module  $\rT_3(E)$ is a free module of rank $1$.
Therefore we obtain a character 
\[
\chi_E : \Gal_{\bQ(\zeta_3)} \to R^\ast
\]
describing the action on $\rT_3(E)$, and a posteriori also the action on $E_{3^n}$. Since $E/\bQ$ has good reduction outside $v \mid 30$, the character $\chi_E$ is unramified outside $2,3,5$.
\begin{lem} \label{lem:1}
The character $\chi_E$ is surjective onto the $1$-units $U^1_R = 1+(\zeta_3 -1)R \subset R^\ast$.
\end{lem}
\begin{pro}
First, the action is via a $3$-group since $E_{3^n}$ has a filtration with quotients 
\[
E_{3^{m+1}}/E_{3^m} \cong E_3,
\]
which over $\bQ(\zeta_3)$ is further filtered with trivial quotients by \eqref{eq:3torsofjacofselmer}. So the action is unipotent and therefore through a $3$-group. Hence $\chi_E$ takes values in $U^1_R$. The group of $1$-units has the structure
\[
U^1_R \cong \mu_3 \times \bZ_3 \times \bZ_3.
\]
By Nakayama's Lemma it suffices to check the surjectivity of  
\[
\im(\chi_E) \surj U^1_R/(U^1_R)^3.
\]
A computation reveals that the image of $U^1_R$ in $(R/9R)^\ast$
is killed by $3$: we compute modulo $9$ and use $3= (\zeta_3+1)(\zeta_3-1)^2$ to see 
\[
\big(1+(\zeta_3-1)a\big)^3 = 1 + 3(\zeta_3-1)a(1+(\zeta_3-1)a) + (\zeta_3-1)^3 a^3
\]
\[
= 1 + (\zeta_3-1)^3\big( (\zeta_3+1)a(1+(\zeta_3-1)a) + a^3\big)
\]
\[
\equiv 1 + (\zeta_3 -1)^3(2a + a) \equiv 1.
\]
The induced map  
\[
U^1_R/(U^1_R)^3 \to (R/9R)^\ast
\]
is surjective and both sides have $27$ elements, hence an isomorphism. It therefore remains to show that $\Gal_{\bQ(\zeta_3)}$ acts via a group of order $27$ on $E_9$. A computation with SAGE \cite{sage} shows
\[
\bQ(E_3) = \bQ(\zeta_3, \sqrt[3]{60}),
\]
\begin{equation} \label{eq:jacofselmer9tors}
\bQ(E_9) = \bQ(\zeta_9, \sqrt[3]{60}, \sqrt[3]{3}).
\end{equation}
and this completes the proof.
\end{pro}

\begin{lem} \label{lem:3}
Any non-trivial $R$-submodule of $E_{3^n}$ contains the $\mu_3$-part of $E_3$.
\end{lem}
\begin{pro}
As an $R$-module $E_{3^n}$ is cyclic, hence every submodule contains the $(\zeta_3-1)$-torsion, which is a module of $3$-elements contained in the $3$-torsion $E_3$. By the Lemma 1 above, an $R$-submodule is also a $\Gal_{\bQ(\zeta_3)}$-submodule, and there is only $\mu_3 \subset E_3$, becasue $60$ is not a cube in $\bQ(\zeta_3)$.
\end{pro}

\begin{cor} \label{cor:4}
Let $v$ be a place different from $2,3,5$ that splits in $\bQ(\zeta_3)$. If $\mu_3 \subseteq E_3$ is not in the kernel of the map 
\[
E_3/(\Frob_v -1)E_3 \to E_{3^n}/(\Frob_v-1)E_{3^n},
\]
induced by inclusion $E_3 \subset E_{3^n}$, 
then $\Frob_v =1$ in $\Gal(\bQ(E_{3^n})/\bQ)$.
\end{cor}
\begin{pro}
As we assume that $v$ splits in $\bQ(\zeta_3)$, we know that $\Frob_v$ acts via $R$-linear automorphisms. In particular the subgroup 
\[
(\Frob_v-1)E_{3^n} \subset E_{3^n}
\]
is an $R$-submodule. Lemma  \ref{lem:3} shows that either $\mu_3 \subseteq (\Frob_v-1)E_{3^n}$, or $(\Frob_v-1)E_{3^n} = 0$.
\end{pro}

\medskip

We recall from \cite{ciperianistix:bashmakov} \S8 the result of the computation of the following \'etale cohomology groups (see loc.\ cit.\ for the definition of $\rH^i_!(U,E_3)$)
\[
\rH^1(U,E_3)  = \rH^1_{\Sel_3}(\bQ,E_3) = \rH^1_{\Sel}(\bQ,E_3) = \Sha(E/\bQ) \cong \bZ/3\bZ \times \bZ/3\bZ,
\]
\[
\rH^1_{\rc}(U,E_3)  = \rH^1_!(U,E_3) = \rH^1_{\Sel^3}(\bQ,E_3) =  \langle [S] \rangle \cong \bZ/3\bZ,
\]
and
\[
\rH^1(U,E_{3^n})  = \rH^1_{\Sel_3}(\bQ,E_{3^n}) \cong \bZ/3\bZ \times \bZ/3^n\bZ,
\]
\[
\rH^1_{\rc}(U,E_{3^n})  = \rH^1_!(U,E_{3^n}) = \rH^1_{\Sel^3}(\bQ,E_{3^n}) =  \langle [S] \rangle \cong \bZ/3\bZ.
\]
and
\[
\rH^1_{\divisor}(\bQ,E)_{3^\infty} =  \rH^1_{\Sel_3}(\bQ,E_{3^\infty}) = \bZ/3\bZ \times \bQ_3/\bZ_3.
\]
We now compute the obstruction to divisibility by $3^n$ as classes over $U$ by exploiting the cohomology sequence of \'etale cohomology associated to 
\[
0 \to E_{3^n} \to E_{3^{n+1}} \to E_3 \to 0
\]
namely 
\[
0 \to \rH^1(U,E_{3^n}) \to \rH^1(U,E_{3^{n+1}}) \to \rH^1(U,E_3) \xrightarrow{\delta} \rH^2_!(U,E_{3^n}) \to 0.
\]
A priori the map on the right lands in $\rH^2(U,E_{3^n})$. But as we start with $\rH^1(U,E_3)$ that contains by coincidence only locally divisible classes, we obtain a map into $\rH^2_!(U,E_{3^n})$. Then we count and find that the image of $\delta$ is of order $3$, while Artin--Verdier duality says
\[
\# \rH^2_!(U,E_{3^n}) = \# \rH^1_!(U,E_{3^n}) = 3.
\]

\begin{pro}[Proof of  Proposition~\ref{prop:casselsforjacofselmer}]
We only have to deal with "the" $\bZ/3Z$-factor of $\rH^1_{\divisor}(\bQ,E)_{3^\infty}$ and show that it becomes $3^n$ divisible in $\rH^1(V,E)$ for small enough $V \subset U$. This factor is generated by an element in $\rH^1(U,E_3)$ and as such has an obstruction in $\rH^2_!(U,E_{3^n})$ which by the above is non-zero. Well, as an obstruction against being $3^n$-divisible in $\rH^1(U,E)$. So we have to shrink now $U$ to $V$. Then the obstruction lies in the dual of the subgroup defined by
\[
0 \to \rH^1_!(V,E_{3^n}) \to \rH^1_!(U,E_{3^n}) \to \bigoplus_{v \in U \setminus V} \rH^1_{\nr}(k_v,E_{3^n}).
\]
So we have to find a $V \subset U$ such that  the Selmer curve that generates $\rH^1_!(U,E_{3^n})$ does not die in 
\[
 \rH^1_{\nr}(k_v,E_{3^n}) = E_{3^n}/(\Frob_v-1)E_{3^n}
\]
upon evaluating a representing cocycle in the Frobenius $\Frob_v$.

Here we have to recall, where the cocycles take their values: we have a surjection
\[
\langle 2,3,5 \rangle = \bZ[1/30]^\ast/(\bZ[1/30]^\ast)^3 = \rH^1(U,\mu_3) \surj \rH^1(U,E_3)
\]
and the Selmer curve is represented by the class of $6$, it means that evaluation in Frobenius takes values in the $\mu_3 \subset E_3 \subset E_{3^n}$. Now, a first conditon for non-trivial value is that $\mu_3$ does survive in 
\[
E_3/(\Frob_v -1)E_3.
\]
But if $\Frob_v$ acts non-trivially on $E_3$, then the coinvariants are a proper quotient, namely the $\bZ/3\bZ$-quotient, and therefore $\mu_3$ dies. Thus we must have that $v$ is chosen in such a way that $v$ splits completely in $\bQ(E_3)/\bQ$. Then the Selmer curve localises to 
\[
\text{Selmer curve cocycle}(\Frob_v) = \Frob_v(\sqrt[3]{6})/\sqrt[3]{6} \in \mu_3 \subset E_3.
\]
We further want that this survives under mapping to $E_{3^n}/(\Frob_v-1)E_{3^n}$ and by Lemma \ref{cor:4} this means that we need $\Frob_v=1$ in $\Gal(\bQ(E_{3^n})/\bQ)$.

We see that we can find $V = U \setminus \{v\}$ with $\rH^1_!(V,E_{3^n}) = 0$ if and only if 
\[
\bQ(E_{3^n}) \text{ and } \bQ(\sqrt[3]{6},\zeta_3)
\]
are linearly disjoint over $\bQ(\zeta_3)$. Now over $\bQ(\zeta_3)$ both fields are abelian and it therefore suffices to check linear disjointness with the field extension in $\bQ(E_{3^n})/\bQ(\zeta_3)$ that correspoinds to the mod $3$ quotient of the Galois group, that is $\bQ(E_9)/\bQ(\zeta_3)$ as follows from the proof of Lemma~\ref{lem:1}. But in \eqref{eq:jacofselmer9tors} we have computed  $\bQ(E_9)$ explicitly, and it is not difficult to see that 
\[
\sqrt[3]{6} \notin \bQ(E_9).
\]
This concludes the proof of Proposition~\ref{prop:casselsforjacofselmer}.
\end{pro}

%%%%%%%%%%%%%%%%%%%%%%%%%%%%%%%%%%%%%%%%%%
%%%%%% Bibliography %%%%%%%%%%%%%%%%%%%%%%%%%%%%%
%%%%%%%%%%%%%%%%%%%%%%%%%%%%%%%%%%%%%%%%%%

\end{document}